\ifdefined\pdfoutput
\pdfoutput=1
\fi
\documentclass[a4paper,11pt]{article}
\usepackage[margin=2.5cm]{geometry}
\usepackage{etoolbox}
\usepackage{xparse}
\usepackage[utf8]{inputenc}
\usepackage[T1]{fontenc}
\usepackage{amsmath}
\usepackage{amssymb,amsfonts,amsopn,dsfont,amscd}
\usepackage{amsthm}
\usepackage{latexsym}
\usepackage{bm}
\usepackage{mathtools}
\usepackage{nicefrac}
\usepackage{mathrsfs}
\usepackage{cancel}
\usepackage{url}
\usepackage{stmaryrd}
\usepackage{graphicx}
\usepackage{natbib}
\usepackage{authblk}
\usepackage[hidelinks]{hyperref}
\usepackage[noabbrev, capitalize]{cleveref}

\newcommand{\C}{\mathcal{C}}
\newcommand{\bC}{\mathbf{C}}

\newcommand{\R}{\mathbb{R}}
\newcommand{\Z}{\mathbf{Z}}
\newcommand{\X}{\mathbf{X}}
\newcommand{\Curl}{\nabla \times }
\newcommand{\Grad}{\nabla}
\newcommand{\Div}{\nabla \cdot}
\newcommand{\ttrace}{\boldsymbol{\gamma}_{\parallel}}

\newcommand{\lift}{L_h}
\newcommand{\elproj}{\Pi_h}

\newcommand{\Tmap}{T_h}
\newcommand{\hTmap}{\widehat{T}_h}

\newcommand{\VectorLtwo}{\mathbf{L}^2(\Omega)}
\newcommand{\VectorLtwoGamma}{\mathbf{L}^2(\Gamma)}
\newcommand{\VectorLtwoGammaPar}{\mathbf{L}^2_{\parallel}(\Gamma)}
\newcommand{\Hcurl}{\mathbf{H}(\operatorname{curl}, \Omega)}
\newcommand{\Hzcurl}{\mathbf{H}_0(\operatorname{curl}, \Omega)}
\newcommand{\Hzcurlz}{\mathbf{H}_0(\operatorname{curl}_0, \Omega)}
\newcommand{\Hdiv}{\mathbf{H}(\operatorname{div}, \Omega)}
\newcommand{\Hzdiv}{\mathbf{H}_0(\operatorname{div}, \Omega)}
\newcommand{\Hzdivz}{\mathbf{H}_0(\operatorname{div}_0, \Omega)}

\newcommand{\VectorHs}{\mathbf{H}^s(\Omega)}

\newcommand{\VectorHtwo}{\mathbf{H}^2(\Omega)}
\newcommand{\VectorHsmhPar}{\mathbf{H}^{s-\frac{1}{2}}_{\parallel}(\Gamma)}
\newcommand{\Vh}{\mathbf{V}_h}
\newcommand{\Qh}{Q_h}

\newcommand{\h}{\mathfrak{H}}

\newcommand{\jump}[1]{\llbracket #1 \rrbracket}

\DeclareMathOperator{\tr}{tr}

\newcommand{\bfu}{\mathbf{u}}
\newcommand{\bfv}{\mathbf{v}}
\newcommand{\bfw}{\mathbf{w}}
\newcommand{\bfz}{\mathbf{z}}
\newcommand{\bff}{\mathbf{f}}
\newcommand{\bfg}{\mathbf{g}}

\newcommand{\Hndivz}{\mathbf{H}^1_n(\operatorname{div}_0, \Omega)}

\newcommand{\harm}{\mathfrak{h}}

\newcommand{\bfzero}{\mathbf{0}}
\newcommand{\nvec}{\mathbf{n}}
\newcommand{\symgrad}{\boldsymbol{\epsilon}}

\renewcommand{\epsilon}{\ensuremath\varepsilon}

\renewcommand{\phi}{\ensuremath{\varphi}}

\theoremstyle{plain}
\newtheorem{theorem}{Theorem}
\newtheorem{corollary}[theorem]{Corollary}
\newtheorem{lemma}[theorem]{Lemma}

\theoremstyle{definition}

\theoremstyle{remark}
\newtheorem{remark}{Remark}

\numberwithin{theorem}{section}
\numberwithin{equation}{section}
\numberwithin{remark}{section}

\setcitestyle{square}

\title{H(curl)-based approximation of the Stokes problem with slip boundary conditions}

\author[1]{Wietse M. Boon}
\author[2]{Ralf Hiptmair}
\author[2]{Wouter Tonnon}
\author[3]{Enrico Zampa\thanks{Corresponding author: \href{mailto:enrico.zampa@univie.ac.at}{enrico.zampa@univie.ac.at}}}
\affil[1]{Faculty of Mathematics, University of Duisburg-Essen, 45127 Essen, Germany}
\affil[2]{SAM, ETH Z\"{u}rich, CH-8092, Switzerland}
\affil[3]{Department of Mathematics, University of Vienna, 1090, Austria}
\date{}

\begin{document}

\maketitle

\begin{abstract}
  Reformulating the incompressible Stokes equations with the velocity sought in
  $\mathbf{H}(\text{curl})$ has recently emerged as a promising approach for the design of
  helicity-preserving schemes in magnetohydrodynamics and pressure-robust finite element
  methods on polygonal meshes. A key challenge in this setting, however, is the treatment
  of Navier slip boundary conditions. In this paper, we overcome this difficulty by
  recasting the slip condition as a Robin boundary condition and proving well-posedness
  (up to a finite-dimensional kernel) of the resulting continuous problem. We further identify the geometric and regularity
  assumptions on the domain and the exact solution under which the classical Stokes
  solution is recovered. Finally, we study a
  finite element Galerkin
  discretization, establishing stability and a priori error estimates. Numerical
  experiments validate the convergence rates predicted by the theory.

\end{abstract}

\section{Introduction}
 The use of $\Hcurl$-conforming elements for fluid problems was pioneered by
  \citet{Girault88, Girault90}. In the wake of the wide adoption of finite element
  exterior calculus (FEEC), this unconventional approach has gained
  popularity. The FEEC approach suggests that physical quantities which can be modeled by
  means of differential forms should be approximated by discrete differential
  forms. Thus, if a velocity is viewed as a $1$-form, its natural finite element
  approximation should rely on discrete 1-forms, that is, $\Hcurl$-conforming finite element
  spaces. This principle has been applied to the semi-Lagrangian discretization of
  the Navier-Stokes equations \citep{TonnonHiptmairSemiLagrangianNavierStokes}, dual field
  formulations \citep{DualFieldNS, ShipengRuijie}, pressure-robust polygonal methods
  \citep{BdVDassiDiPietroDroniou22, DiPietroDroniuQian24} and recent advancements in
  magnetohydrodynamics (MHD), where these elements have proven particularly effective for
  preserving cross helicity \citep{HuLeeXu21, LaakmanHuFarrell23}, see also
  \citep{AlonsoRapetti26}.

However, a significant omission in these works is the treatment of general boundary
  conditions. A preliminary step toward addressing this gap was taken in \citep{Nitsche},
  which introduced a Nitsche-type strategy to impose no-slip boundary conditions. In the
  present work, we address Navier slip boundary conditions. For simplicity, we consider
  the stationary Stokes equations; however, the proposed boundary treatment naturally
  extends to more general problems.

We start by formulating the continuous problem and show that it is well-posed
  up to a finite-dimensional kernel using the theory of compact perturbations of coercive operators on the space
  $\X=\Hcurl\cap \nabla H^1(\Omega)^\perp$. We also show that the continuous formulation
  recovers the correct physical solution under suitable regularity assumptions. As our
  method reduces to the aforementioned $\Hcurl$-based formulations on polyhedra, we
  conclude that these formulations may fail to recover the correct Stokes solution on
  non-convex domains.

We then consider the finite element Galerkin discretization, for which we cannot
  directly apply the standard theory, since the discrete space contains elements of
  $\Hcurl$ that are orthogonal only to discrete gradients. This means that the discrete
  space is non-conforming, as it is not a subspace of $\X$. Fortunately, the discrete and
  continuous spaces can be shown to be \lq\lq close enough\rq\rq\, to prove well-posedness and
  convergence.

\textbf{Outline.} In \Cref{sec:VVP} we recall how to derive the
  rotational formulation of the Stokes problem and an equivalent expression for Navier slip
  boundary conditions involving the Weingarten map due to \citet{MitreaMonniaux}. In
  \Cref{sec:vf} we write the continuous variational formulation and show its
  well-posedness. In \Cref{sec:discrete_vf} we focus on the discrete formulation, derive a
  priori error estimates for its solution, and extend the analytical results to the
  saddle-point system involving the pressure. Finally, in \Cref{sec:numerics}, we describe
  the implementation of the method and validate it through a series of two- and
  three-dimensional numerical tests. In particular, we demonstrate that the method
  performs well when using either the exact Weingarten map or an approximate one computed
  as in \citep{NeunteufelCurvatureComputation}, even on non-curved meshes. This provides a
  novel and elegant approach to overcoming Babu\v{s}ka's paradox in the context of slip
  boundary conditions; see \citep{Verfurth}.
\section{The Stokes problem in rotational form with
  Navier slip boundary conditions}

\label{sec:VVP}
Let $\Omega\subset \R^d$, $d = 2, 3$, be a domain and denote by $\Gamma$ its boundary. Formally, we seek a velocity field $\bfu: \Omega\to\R^d$ and pressure field $p: \Omega\to\R$ such that
\begin{subequations}\label{eq:continuousSymGradStokes}
\begin{align}
    - 2 \Div\symgrad(\bfu) + \Grad p &= \mathbf{f}, & \text{ in }&\Omega, \label{eq:continuousSymGradStokesa} \\
    \Div \bfu &= 0,  & \text{ in }&\Omega, \\
    \label{eq:unBC}
    \bfu\cdot \nvec &= 0, & \text{ on }&\Gamma,\\
    [\symgrad(\bfu)\nvec]_{\mathrm{tan}} &= 0, & \text{ on }&\Gamma,
    \label{eq:NavierSlipBC}
    \end{align}
\end{subequations}
where $\bff:\Omega\to\R^d$ is a given forcing term and $\nvec:\Gamma\to\R^d$ is
the outward oriented unit vector normal to $\Gamma$. Moreover, $\symgrad(\bfu)$ denotes
the symmetric gradient of $\bfu$, that is,
\begin{equation*}
    \symgrad(\bfu) \coloneqq \frac12 \left(\nabla \bfu + (\nabla\bfu)^T \right).
\end{equation*}
Finally,
$\boldsymbol{\zeta}_{\mathrm{tan}}\coloneqq \boldsymbol{\zeta} - (\boldsymbol{\zeta}\cdot
\nvec)\nvec$ stands for the tangential part of a vector field $\boldsymbol{\zeta}$ on
$\Gamma$. Note that the boundary condition in \cref{eq:NavierSlipBC} arises as a natural
  boundary condition from integration by parts of the diffusion term in
\cref{eq:continuousSymGradStokesa}.
After fixing the additive constant in $p$, the solution to \eqref{eq:continuousSymGradStokes} is
unique, unless $\Omega$ is axisymmetric. In that case, we face a non-trivial kernel
consisting of rigid body rotations about the axis of symmetry. Precise statements, including the compatibility condition on the load and the case of balls and spherical shells, are given in \Cref{thm:FiniteDimensionalKernel,thm:kernelStrainForm}.

In the case that $\Gamma$ is a $\C^2$-boundary, \citet[Eq. (2.9)]{MitreaMonniaux} have
shown that the homogeneous Navier slip boundary conditions \cref{eq:unBC} and
\cref{eq:NavierSlipBC} are equivalent to
\begin{align} \label{eq: MitreaMonniaux}
    \bfu \cdot \nvec &= 0, &
    -\nvec \times (\Curl \bfu) +2\mathcal{W}(\bfu_{\mathrm{tan}}) &= 0,
\end{align}
with $\mathcal{W}:T\Gamma\to T\Gamma$ the Weingarten map, defined by $\mathcal{W}\boldsymbol{\tau}=-D_{\boldsymbol{\tau}}\nvec$, where $D_{\boldsymbol{\tau}}$ denotes the derivative along $\boldsymbol{\tau}$ on $\Gamma$, and $T\Gamma$ the space of tangential vector fields on $\Gamma$.

In the two-dimensional case ($d=2$), the Weingarten map reduces to multiplication by
  the scalar curvature of $\Gamma$. Moreover, in 2D the curl operator is defined as
$\Curl \bfu \coloneqq \partial_x u_y - \partial_y u_x$ for vector-valued
$\bfu = [u_x, u_y]^T$ and $\Curl q \coloneqq [\partial_y q, -\partial_x q]^T$ for scalar
$q$. In two dimensions, planar vectors are identified with vectors in $\R^3$ with vanishing third component and a scalar $q$ with $q\,\mathbf{e}_3$, and cross products are understood in $\R^3$. In particular, $\mathbf{a}\times\mathbf{b}=a_1b_2-a_2b_1$ is a scalar for planar vectors $\mathbf{a},\mathbf{b}$, and $\nvec\times q=q\,[n_2,-n_1]^T$.

For sufficiently smooth $\bfu$ with $\Div\bfu=0$, we can rewrite the second order term in
\Cref{eq:continuousSymGradStokesa} as
\begin{equation} \label{eq: rewrite symgrad}
    \begin{aligned}
        -2 \Div\symgrad(\bfu) &= -\Delta\bfu - \Grad (\Div\bfu)\\
        &= \Curl\Curl\bfu - 2\Grad(\Div\bfu)\\
        &= \Curl\Curl\bfu.
    \end{aligned}
\end{equation}
We conclude formally that the solution of the Stokes problem
\eqref{eq:continuousSymGradStokes} solves the following boundary value problem: find a
velocity field $\bfu: \Omega \to \R^d$ and pressure field $p: \Omega \to \R$ such that
\begin{subequations} \label{eq:strong_problem}
    \begin{align}
        \Curl \Curl \bfu + \Grad p &= \mathbf{f}, & \text{ in }&\Omega, \label{eq:sp_a} \\
        \Div \bfu &= 0,  & \text{ in }&\Omega, \label{eq:sp_b} \\
        \bfu \cdot \nvec &= 0, & \text{ on }&\Gamma, \label{eq:sp_c}  \\
        -\nvec \times (\Curl \bfu) +\alpha(\bfu_{\mathrm{tan}}) &= 0,  & \text{ on }&\Gamma, \label{eq:sp_d}
    \end{align}
\end{subequations}
where $\alpha=2\mathcal{W}$.
Problem \eqref{eq:strong_problem} will be the starting point for our discretization. Note that when $\alpha = 0$, these boundary conditions coincide with those considered e.g. in \citet{TonnonHiptmairSemiLagrangianNavierStokes,BdVDassiDiPietroDroniou22, DiPietroDroniuQian24, ShipengRuijie}, but in general they are different.

  \begin{remark}
    When $\Gamma$ is Lipschitz and \emph{piecewise} of class $\mathcal{C}^2$, the Weingarten map
    is defined on its $C^{2}$-smooth faces. For a velocity field sufficiently smooth up to each
    such face and satisfying $\bfu\cdot\nvec=0$, the identity relating
    \eqref{eq:NavierSlipBC} and \eqref{eq:sp_d} holds pointwise on that face. The
    corresponding integration-by-parts identity on a piecewise $\mathcal{C}^2$ boundary is
    given in \citet[Theorem~3.1.1.2]{Grisvard} for $\mathbf{H}^2$ vector fields. The
    $\mathbf{H}^1$ identity used below is stated in \citet[Theorem~2.3]{CostabelDauge99}.
\end{remark}

\begin{remark} \label{rem:negativityAlpha}
    Note that the determinant of the Weingarten map $\mathcal{W}$ is the Gaussian curvature of $\Gamma$. We can therefore not draw conclusions on the sign of $\alpha$ for general $\Omega$. In the case of a convex three-dimensional $\C^{1,1}$-domain, for example, $\mathcal{W}$ will be negative semi-definite almost everywhere. In general, however, $\mathcal{W}$ may be positive-definite, negative-definite, or indefinite on various parts of the boundary.
\end{remark}

\begin{remark}
  In terms of differential forms, $u$ can be considered a $1$-form or a
  $(d-1)$-form. Equation \eqref{eq:sp_d} then reads either as
  $(\mathrm{d} u)^{\mathrm{nor}} + \alpha u^{\mathrm{tan}}  = 0$ or
  $-(\delta u)^{\mathrm{tan}} + \beta u^{\mathrm{nor}} = 0$, respectively, where $\beta$ denotes $\alpha$ conjugated by the Hodge star of $\Gamma$ (in three dimensions, by the rotation by $90^\circ$ in the tangent plane), so that $\beta=\alpha$ if $\alpha$ is a multiple of the identity. Here
  $\zeta^{\mathrm{tan}}$ and $\zeta^{\mathrm{nor}}$ denote the tangential and the normal
  trace of a differential form $\zeta$,  respectively, see \citep[Def.~2.1 and Table~1]{AFGS}. We then
  immediately see that equation \eqref{eq:sp_d} is a \emph{Robin} boundary condition.
  For the Hodge--Laplace operator for 1-forms, \citet[Sects.~2.5.1--2.5.2]{YangPHD}
    studies this type of Robin boundary condition in the framework of finite element
    exterior calculus, with a positive constant $\alpha$ and the condition
    $\operatorname{tr}(\delta u)=0$ in place of $\bfu\cdot\nvec=0$. The scalar variable
    $\delta u$, denoted $\sigma$ in  that paper, plays the role of a pressure variable in the
    corresponding vector equation. That analysis does not apply to the present problem.
\end{remark}

\begin{remark}
    Problem \eqref{eq:strong_problem} is similar to the time-harmonic Maxwell equations with impedance boundary conditions \citep{GaMe12}. The main difference is the divergence-free constraint which gives additional smoothness to the space in which $\bfu$ is sought, but makes the analysis of the discrete problem more involved. These two issues will be discussed in \Cref{sec:vf,sec:discrete_vf}, respectively.
  \end{remark}

\section{Variational formulation of the continuous problem}
\label{sec:vf}

In this section, we first describe the functional setting and introduce the notation conventions used throughout this work. Afterward, we present and analyze the variational formulation of \eqref{eq:strong_problem}.

\subsection{Functional setting}
Let $\Omega\subset \R^d$, $d = 2, 3$, be a bounded and connected Lipschitz domain. We will sometimes further specify $\Omega$ as either a polyhedron or a curved polyhedron in the sense of \citep{CostabelDauge99}. In any case, we assume that $\Gamma$ is piecewise of class $\mathcal{C}^2$, consisting of finitely many smooth faces $\Gamma_i$, $i = 1, \dots, M$. We denote the $L^2(\Omega)$ inner product by parentheses and the analogous product on its boundary $\Gamma$ by angled brackets:
\begin{align*}
    (\bfu, \bfv) &\coloneqq \int_\Omega \bfu \cdot \bfv\, \mathrm{d} \mathbf{x}, &
    \langle \bfu, \bfv \rangle &\coloneqq \int_\Gamma \bfu \cdot \bfv\,\mathrm{dS}.
\end{align*}
With a slight abuse of notation, we denote inner products for scalar functions using the same parentheses and brackets.

We recall the definitions of the following classical Hilbert spaces
\begin{align*}
    \Hcurl &\coloneqq \{ \bfv \in \VectorLtwo \mid \Curl \bfv \in (L^2(\Omega))^{2d-3}\}, \\
    \Hdiv &\coloneqq \{ \bfv \in \VectorLtwo \mid \Div \bfv \in L^2(\Omega)\},
\end{align*}
and their associated norms
\begin{align*}
    \lVert\bfv\rVert^{2}_{\Hcurl} &\coloneqq \lVert\bfv\rVert^{{2}}_{\VectorLtwo} + \lVert\nabla\times\bfv\rVert^{{2}}_{L^2(\Omega)^{2d-3}},\\
    \lVert\bfv\rVert^{{2}}_{\Hdiv} &\coloneqq \lVert\bfv\rVert^{{2}}_{\VectorLtwo} + \lVert\nabla\cdot\bfv\rVert^{{2}}_{L^2(\Omega)}.
\end{align*}
{We continue by defining the subspaces  $\Hzcurl$ and $\Hzdiv$ as the closure of $\bC_c^{\infty}(\Omega)$ in the $\lVert \cdot \rVert_{\Hcurl}$ and $\lVert \cdot \rVert_{\Hdiv}$ norms, respectively. Moreover,  we define $\Hzcurlz$ and $\Hzdivz$ by
\begin{align*}
    \Hzcurlz &\coloneqq \{ \bfv \in  \Hzcurl\mid \Curl \bfv = 0\}, \\
    \Hzdivz &\coloneqq \{ \bfv \in  \Hzdiv\mid \nabla\cdot\bfv = 0\}.
\end{align*} }

For function spaces with regularity of fractional order, we will use the Sobolev-Slobodeckij seminorm for $0 < s < 1$:
\begin{equation*}
    \lvert \bfv \rvert_{\VectorHs} \coloneqq \left( \int_\Omega\int_\Omega \frac{\lvert \bfv(\mathbf{x})-\bfv(\mathbf{y})\rvert^2}{\lVert\mathbf{x}-\mathbf{y}\rVert^{2s+d}}\,d\mathbf{x}\,d\mathbf{y}\right)^{{\frac{1}{2}}}.
\end{equation*}
This allows us to define the space
\begin{align*}
    \VectorHs &\coloneqq \{\bfv\in\VectorLtwo\mid  \lvert \bfv \rvert_{\VectorHs}<\infty\},
\end{align*}
for $0<s<1$. For non-integer $s > 1$, let $\sigma = \lfloor s\rfloor$. Then $\VectorHs$ is the space of functions in $\mathbf{H}^\sigma(\Omega)$ whose weak derivatives of order $\sigma$ are in $\mathbf{H}^{s - \sigma}(\Omega)$. For integer $s$, $\VectorHs$ is defined as usual through weak derivatives, with $\mathbf{H}^0=\mathbf{L}^2$.
{We now define some function spaces on $\Gamma$, following \citep{BuffaCostabelSheen}. We start with}
\begin{equation*}
    \VectorLtwoGammaPar\coloneqq \{ \boldsymbol{\zeta}\in \VectorLtwoGamma\mid \boldsymbol{\zeta}\cdot \nvec = 0\}.
\end{equation*}
We are now ready to define the tangential {component} $\ttrace$ and the trace operator $ \boldsymbol{\gamma}_t$: $\ttrace: \bC^{0}(\overline{\Omega})\to\VectorLtwoGammaPar$ as
\begin{align*}
    \ttrace(\bfv) &\coloneqq
    \bfv|_{\Gamma} - (\bfv|_{\Gamma} \cdot \nvec) \nvec, \qquad \boldsymbol{\gamma}_t(\bfv)\coloneqq \bfv|_{\Gamma}\times \nvec, \qquad&
    \forall \bfv &\in \bC^{0}(\overline{\Omega}).
\end{align*}
{We need $\boldsymbol{\gamma}_t$ only to define the range of $\ttrace$. As stated in \cite[Thm.~3.10]{EG1}, $\ttrace$ and $ \boldsymbol{\gamma}_t$ can be extended to bounded maps on $\VectorHs$, $\frac{1}{2} < s \leq \frac{3}{2}$. In particular, we define
\begin{equation*}
    \mathbf{H}^{s-\frac{1}{2}}_{\parallel}(\Gamma) \coloneqq \ttrace(\mathbf{H}^s(\Omega)), \qquad \mathbf{H}_t^{s-\frac{1}{2}}(\Gamma) \coloneqq \boldsymbol{\gamma}_t(\mathbf{H}^s(\Omega)).
\end{equation*}
These spaces are equipped with the associated quotient norms, see \citep[Def.~2.2]{BuffaCostabelSheen}. We denote by $\mathbf{H}_t^{\frac{1}{2}-s}(\Gamma)$ the dual space of $\mathbf{H}_t^{s-\frac{1}{2}}(\Gamma)$. Then, it was shown by \citet{BuffaCostabelSheen} that $\ttrace$ can be extended to a bounded map $\ttrace:\Hcurl\to \mathbf{H}^{-\frac{1}{2}}_{t}(\Gamma)$. We will make a slight abuse of notation and not distinguish between these two extensions of $\ttrace$.} We summarize these results in the following lemma.
Finally, if the faces $\Gamma_i$ are of class $\C^{r+1}$ for an integer $r\geq1$, we denote by $\mathbf{H}^r_{\mathrm{pw}}(\Gamma)$ the space of fields $\boldsymbol{\zeta}\in\VectorLtwoGammaPar$ with $\lVert\boldsymbol{\zeta}\rVert_{\mathbf{H}^r_{\mathrm{pw}}(\Gamma)}<\infty$, where
\begin{equation*}
    \lVert\boldsymbol{\zeta}\rVert^2_{\mathbf{H}^r_{\mathrm{pw}}(\Gamma)}\coloneqq\sum_{i=1}^M\lVert D\boldsymbol{\psi}_i^T(\boldsymbol{\zeta}\circ\boldsymbol{\psi}_i)\rVert^2_{\mathbf{H}^r(\omega_i)}
\end{equation*}
and $\boldsymbol{\psi}_i\in\mathbf{C}^{r+1}(\overline{\omega_i})$, $\omega_i\subset\R^{d-1}$, are fixed parametrizations of the faces $\Gamma_i$ with $D\boldsymbol{\psi}_i$ of full rank. This is a broken norm: no regularity of $\boldsymbol{\zeta}$ across the interfaces between the faces is required.

\begin{lemma} \label{thm:TraceBounded}
    $\ttrace:\VectorHs\to \VectorHsmhPar$, $\frac{1}{2} < s \leq \frac{3}{2}$, and $\ttrace:\Hcurl\to \mathbf{H}^{-\frac{1}{2}}_{t}(\Gamma)$ are bounded.
\end{lemma}

The final space we introduce is
\begin{gather*}
 \X \coloneqq  \Hcurl\cap \Hzdivz,
\end{gather*}
which we equip with the $\Hcurl$-norm. Note that $\X$ can equivalently be characterized as
\begin{equation} \label{eq: def X}
    \X = \{ \bfv \in \Hcurl\mid (\bfv, \Grad q) = 0, \forall q\in H^1(\Omega)\}.
\end{equation}
{Let $\h^1\coloneqq \{ \bfv\in \X\mid \Curl \bfv = 0\}$ be the space of Neumann harmonic vector fields and let $\Z$ be its orthogonal complement in $\X$. This immediately yields the Hodge-Helmholtz decomposition,
\begin{equation*}
    \X = \h^1 \oplus^{\perp} \Z,
\end{equation*}
because $\Hcurl = \Grad H^1(\Omega) \oplus^{\perp} \h^1 \oplus^{\perp} \Z$.
}
Recall that the Poincar\'e-{Steklov} inequality holds on {$\Z$}:
\begin{theorem}
\label{thm:PS}
There exists a constant $C_P>0$ such that
\begin{equation*}
  \lVert \bfz\rVert_{\VectorLtwo}\leq C_P\lVert \Curl \bfz\rVert_{L^2(\Omega)^{2d-3}} \quad {\forall \bfz\in \Z}.
\end{equation*}
\end{theorem}
\begin{proof}
    {The result follows from the compact embedding of $\X$ in $\VectorLtwo$ \citep{Weber}.  See also \citep[Corollary~4.4]{Hip02}}.
\end{proof}
We proceed with some regularity properties.
\begin{theorem} \label{thm:embeddingXinHs} If $\Omega$ is a Lipschitz polyhedron, there
  exists $s>\frac{1}{2}$ such that $\X$ is continuously embedded in $\VectorHs$. If
  $\Omega$ is either a convex domain or has a $\C^{1,1}$ boundary, this statement holds
  with $s = 1$.
\end{theorem}
\begin{proof}
See \citep[Thm. 2.9, 2.17, Prop. 3.7]{Amr+98}.
\end{proof}

From now on, {we will assume that $\Omega$ is such that $\X$ is embedded in $\mathbf{H}^s(\Omega)$ with $s > \frac{1}{2}$.}
As a consequence of {this assumption,} we obtain the following trace inequality in $\X$.

\begin{corollary} \label{cor:trace_ineq}
There exists $C_{\tr}>0$ such that for all $\bfv\in \X$
\begin{equation*}
    \lVert \ttrace(\bfv)\rVert_{\VectorLtwoGamma}\leq C_{\tr} { \lVert \bfv\rVert_{\Hcurl}}.
\end{equation*}
\end{corollary}
\begin{proof}
   Let $s \in (1/2, 1]$ be
the regularity index in the standing assumption.  We have for $\bfv\in \X$:
    \begin{align*}
        \lVert \ttrace(\bfv)\rVert_{\VectorLtwoGamma} \leq C\lVert \ttrace(\bfv)\rVert_{\VectorHsmhPar}
        \leq C\lVert \bfv\rVert_{\VectorHs}
        \leq C{ \lVert \bfv\rVert_{\Hcurl}},
    \end{align*}
    where we used \Cref{thm:TraceBounded} for the second inequality and the standing assumption for the third inequality.
\end{proof}

{We conclude this section with an auxiliary result regarding the space $\h^1$.
\begin{theorem}
\label{thm:harm}
The quantities $\lVert \cdot \rVert_{\mathbf{H}^s(\Omega)}$, $\lVert \ttrace( \cdot )\rVert_{\mathbf{L}^2(\Gamma)}$, and $\lVert \cdot \rVert_{\Hcurl}$ are all equivalent norms on $\h^1$.
\end{theorem}
\begin{proof}
We start by showing that $\lVert\ttrace(\cdot)\rVert_{\mathbf{L}^2(\Gamma)}$ is a norm on
$\h^1$. Let $\bfv\in \h^1$ and assume $\ttrace(\bfv) = 0$. Then $\bfv$ belongs to
$\mathbf{H}_0(\operatorname{curl}_0, \Omega) \cap \Hzdivz$, which is a subspace of $\mathbf{H}^1_0$ by \citep[Theorem 2.5]{Amr+98}. Then, \citep[Remark 2.6, Eq.~(2.9)]{Amr+98} implies that $\bfv = 0$, proving the claim. The proof is concluded by invoking equivalence of all norms in finite dimensional spaces.
\end{proof}}

\subsection{The abstract continuous problem}

Let $\alpha$ be a field on $\Gamma$ such that its restriction to a smooth face $\Gamma_i$ is an essentially bounded, self-adjoint endomorphism of the tangent bundle of $\Gamma_i$, for $i = 1, \dots, M$. We denote by $\lVert\alpha\rVert_{L^\infty(\Gamma)}$ the essential supremum of its operator norm.  Multiplying equation \eqref{eq:sp_a} by $\bfv$, integrating by
parts and using the boundary condition \eqref{eq:sp_d} we obtain the continuous problem: given $\bff \in \VectorLtwo$,
find $\bfu\in \X$ such that
\begin{equation}
    a(\bfu, \bfv) = (\mathbf{f}, \mathbf{v})\quad
    \forall \bfv \in \X.
    \label{eq:continuous_problem}
\end{equation}
Here we have defined the bilinear forms $a, \widetilde{a}, k: \X\times \X\to \R$ as
\begin{gather}
  \label{eq: bilinear forms}
  \begin{aligned}
    a &\coloneqq \widetilde{a} + k,&\widetilde{a}(\bfu, \bfv) &\coloneqq (\Curl \bfu,
    \Curl \bfv) {+ \langle \ttrace(\bfu), \ttrace(\bfv) \rangle}, \\ && k(\bfu, \bfv)
    &\coloneqq \langle {(\alpha - I)} \ttrace(\bfu), \ttrace(\bfv)\rangle.
  \end{aligned}
\end{gather}

We analyze \eqref{eq:continuous_problem} showing several key properties of the bilinear forms in the following lemma.
\begin{lemma} \label{lem:boundednessAndCoercivityBilinearForms}
The bilinear forms $\widetilde{a}, k$ and $a$ are symmetric and bounded on $\X$. Moreover, $\widetilde{a}$ is coercive on $\X$, that is, there exists $\gamma>0$ such that for all $\bfv\in\X$
\begin{equation*}
    \widetilde{a}(\bfv,\bfv) \geq\gamma\lVert \bfv\rVert_{{\Hcurl}}^2.
\end{equation*}
\end{lemma}
\begin{proof}
    The symmetry of the bilinear forms is immediate due to the assumed symmetry {of} $\alpha$.
    For boundedness of $k$, we first recall that $\alpha$ is essentially bounded on $\Gamma$. For $\bfu,\bfv\in \X$, \Cref{cor:trace_ineq} implies,
    \begin{align*}
        \lvert k(\bfu,\bfv) \rvert &\leq (\lVert \alpha \rVert_{L^\infty(\Gamma)} { + 1})\lVert  \ttrace(\bfu)\rVert_{\VectorLtwoGamma} \lVert \ttrace(\bfv)\rVert_{\VectorLtwoGamma}\\
        &\leq C (\lVert \alpha \rVert_{L^\infty(\Gamma)} { + 1} ){\lVert \bfu \rVert_{\Hcurl} \lVert \bfv \rVert_{\Hcurl} }.
    \end{align*}
    {The boundedness of $\widetilde{a}$ can be shown similarly.} Thus, we conclude boundedness of
    $a: \X \times \X \to \R$. {We conclude by showing coercivity of $\widetilde{a}: \X\times \X\to\R$. Let $\bfv \in \X$. Then we can write $\bfv = \bfv^{\harm} + \bfv^{\perp}$ with $\bfv^{\harm}\in \h^1$ and $\bfv^{\perp}\in \Z$. Then, combining \Cref{thm:PS} and \Cref{thm:harm} we obtain
    \begin{align*}
        \lVert \bfv \rVert_{\Hcurl}^2 &= \lVert \bfv^{\harm} + \bfv^{\perp} \rVert_{\Hcurl}^2 \\
        & = \lVert \bfv^{\harm} + \bfv^{\perp} \rVert_{\VectorLtwo}^2 + \lVert \Curl\bfv^{\perp} \rVert_{L^2(\Omega)^{2d-3}}^2\\
        &= \lVert \bfv^{\harm} \rVert^2_{\VectorLtwo} + \lVert \bfv^{\perp} \rVert^2_{\VectorLtwo} + \lVert \Curl\bfv^{\perp} \rVert_{L^2(\Omega)^{2d-3}}^2\\
        & \leq C\lVert \ttrace(\bfv^{\harm} ) \rVert_{\mathbf{L}^2(\Gamma)}^2 +( 1 + C_P^2) \lVert \Curl \bfv^{\perp} \rVert_{L^2(\Omega)^{2d-3}}^2 \\
        & \leq C\lVert \ttrace(\bfv)\rVert_{\mathbf{L}^2(\Gamma)}^2  + C \lVert \ttrace(\bfv^{\perp} ) \rVert_{\mathbf{L}^2(\Gamma)}^2 +  ( 1 + C_P^2) \lVert \Curl \bfv^{\perp} \rVert_{L^2(\Omega)^{2d-3}}^2 \\
        & \leq C\lVert \ttrace(\bfv)\rVert_{\mathbf{L}^2(\Gamma)}^2  + CC_{\mathrm{tr}}^2 \lVert \bfv^{\perp} \rVert_{\Hcurl}^2 +  ( 1 + C_P^2) \lVert \Curl \bfv^{\perp} \rVert_{L^2(\Omega)^{2d-3}}^2\\
        & \leq  C\lVert \ttrace(\bfv)\rVert_{\mathbf{L}^2(\Gamma)}^2  + (1 + C_P^2)(1 + CC_{\mathrm{tr}}^2) \lVert \Curl \bfv \rVert_{L^2(\Omega)^{2d-3}}^2\\
        & \leq C \widetilde{a}(\bfv, \bfv).
    \end{align*}
    This proves coercivity.}
\end{proof}

With \Cref{lem:boundednessAndCoercivityBilinearForms}, we prove the following result.
\begin{theorem} \label{thm:FiniteDimensionalKernel}
    The null space
    \begin{equation}
        \ker (a)\coloneqq \{ \bfu \in \X \mid a(\bfu, \bfv) = 0, \forall \bfv \in \X \},
        \label{eq: kernel_of_a}
    \end{equation}
    is finite-dimensional. Moreover, Problem \eqref{eq:continuous_problem} is solvable if and only if
\begin{equation*}
 (\bff,\bfz)=0\qquad\text{for every }\bfz\in\ker(a).
\end{equation*}
When this compatibility condition holds, the solution is unique modulo $\ker(a)$.
The representative orthogonal to $\ker(a)$ in the $\Hcurl$ inner product satisfies
\begin{equation*}
 \lVert\bfu\rVert_{\Hcurl}\le C\lVert\bff\rVert_{\VectorLtwo}.
\end{equation*}
In particular, a trivial kernel gives a unique solution for every
$\bff\in\VectorLtwo$.
\end{theorem}
\begin{proof}
    Define the operators $A, \widetilde{A}, K:\X \to \X^*$ via
    \begin{align*}
        \langle A\bfu, \bfv\rangle_{\ast} &\coloneqq a(\bfu, \bfv), &
        \langle \widetilde{A}\bfu, \bfv\rangle_{\ast} &\coloneqq \widetilde{a}(\bfu, \bfv), &
        \langle K\bfu, \bfv\rangle_{\ast} &\coloneqq k(\bfu, \bfv),
    \end{align*}
Here $\bfu,\bfv\in\X$, and $\X^*$ denotes its dual.
By \Cref{lem:boundednessAndCoercivityBilinearForms} and the Lax--Milgram theorem,
$\widetilde A$ is an isomorphism from $\X$ to $\X^*$.
To prove compactness, choose $t$ with $1/2<t<s$ and factor the trace map as
\begin{equation*}
 T=\ttrace|_{\X}:\quad
 \X\longrightarrow\mathbf H^s(\Omega)
 \lhook\joinrel\longrightarrow\mathbf H^t(\Omega)
 \longrightarrow\VectorLtwoGammaPar.
\end{equation*}
The middle map is compact, and the other two maps are bounded. Thus $T$ is compact;
this is the argument used in \citep[proof of Lemma~4.5]{GaMe12}.
Multiplication by $\alpha-I$ is bounded on $\VectorLtwoGammaPar$, so
$K=T^*M_{\alpha-I}T$ is compact. Hence
\begin{equation*}
 A=\widetilde A(I+\widetilde A^{-1}K)
\end{equation*}
is Fredholm of index zero; see \citep{BrezisBook} for these standard operator facts.
Its kernel is finite-dimensional. Symmetry gives
$\langle A\bfv,\bfz\rangle_*=0$ for every $\bfv\in\X$ if and only if
$\bfz\in\ker(a)$. The Fredholm alternative therefore identifies its range as the
functionals annihilating $\ker(a)$, proving the compatibility condition and
uniqueness modulo the kernel. The restriction of $A$ to the closed orthogonal
complement of its kernel has a bounded inverse onto its range. Since the load
$\bfv\mapsto(\bff,\bfv)$ has dual norm at most $\lVert\bff\rVert_{\VectorLtwo}$,
the stated stability estimate follows.
\end{proof}

\begin{remark}
When $\ker(a)$ is nontrivial, both the compatibility condition and a
normalization of the velocity are needed for a uniquely specified solution.
Assuming $\ker(a)=\{0\}$ removes these restrictions on the load and gives
well-posedness on $\X$.
\end{remark}

\begin{remark}
    If $\alpha$ is positive definite, then {$\ker (a)$} is trivial and Problem \eqref{eq:continuous_problem} is well-posed by Lax-Milgram. However, as noted in \cref{rem:negativityAlpha}, $\alpha$ can be a negative or indefinite operator on parts of the boundary, and thus a {deeper} analysis is required.
\end{remark}

{
\subsection{Relationship of \texorpdfstring{$\Hcurl$}{H(curl)}- and \texorpdfstring{$\mathbf{H}^1(\Omega)$}{H1}-based variational formulations}
\label{sub:relation_with_the_stokes_problem}
In this section we assume that $\alpha$ is any $L^{\infty}$ field coinciding with $2 \mathcal{W}$ on the interior of each smooth face $\Gamma_i$, $i = 1, \dots, M$. We are interested in comparing \eqref{eq:continuous_problem} with the strong formulation \eqref{eq:continuousSymGradStokes} and with the more standard weak formulation: find $\bfu \in \Hndivz\coloneqq \{ \bfv \in \mathbf{H}^1(\Omega)\mid \bfv\cdot \nvec = 0\, \text{ on $\Gamma$}, \Div \bfv = 0\}$ satisfying
\begin{equation}
        a^{\prime}(\bfu,\bfv) \coloneqq 2(\symgrad(\bfu), \symgrad(\bfv)) = (\mathbf{f}, \bfv),
        \label{eq:StokesWF}
\end{equation}
for each $\bfv\in \Hndivz$. Note that $\Hndivz\subset \X$, and the two spaces are equal when $\Omega$ is either convex or has a $\C^{1,1}$ boundary as a consequence of \Cref{thm:embeddingXinHs}.

We start by recalling the following important result about the kernel of \eqref{eq:StokesWF}.
\begin{theorem}\label{thm:kernelStrainForm}
Let $\Omega\subset \R^d$ be a bounded connected domain with $\C^{1,1}$ boundary. If $d = 3$, then $\dim \ker(a^{\prime})\leq 3$ and, up to rotations and translations of $\Omega$, $\ker(a^{\prime})$ is spanned (at most) by
\begin{equation*}
    \mathbf{e}_1\times \mathbf{x}, \qquad \mathbf{e}_2\times \mathbf{x}, \qquad \mathbf{e}_3\times \mathbf{x},
\end{equation*}
where $\mathbf{e}_1$, $\mathbf{e}_2$ and $\mathbf{e}_3$ are the canonical basis vectors of $\mathbb{R}^3$.
In particular, only three cases occur:
\begin{equation*}
    \dim \ker(a^{\prime}) = \begin{cases} 0&\text{if $\Omega$ is not axisymmetric;}\\
    1 &\text{if $\Omega$ is axisymmetric about a single axis;}\\
    3 &\text{if $\Omega$ is a ball or a spherical shell.}
    \end{cases}
\end{equation*}
If $d = 2$, then $\dim \ker(a^{\prime})\leq 1$ and, up to translation of $\Omega$, $\ker(a^{\prime})$ is spanned (at most) by $[-y,x]^T$. In particular, $\dim \ker(a^{\prime}) = 1$ if $\Omega$ is either a disk or an annulus, and $\dim \ker(a^{\prime}) = 0$ otherwise.
\end{theorem}
\begin{proof}
The result is proven in \citep{FallocchiGazzola} when $\Gamma$ is connected. The proof can be immediately extended to the general case.
\end{proof}
\begin{theorem}[Relation between weak formulations]
\label{thm:WvsW}
If $\Omega$ is a curved polyhedron and $\bfu\in \Hndivz$ solves \eqref{eq:continuous_problem}, then $\bfu$ solves \eqref{eq:StokesWF}. Conversely, if $\Omega$ is either a convex curved polyhedron or has $\C^{1,1}$ boundary and $\bfu$ solves \eqref{eq:StokesWF}, then it solves \eqref{eq:continuous_problem}.
\end{theorem}
\begin{proof}
Due to the assumption on $\alpha$, the following auxiliary identity holds when $\Omega$ is a curved polyhedron:
    \begin{equation*}
        2(\symgrad(\bfw), \symgrad(\bfv)) = a(\bfw, \bfv), \qquad \forall \bfw, \bfv \in \Hndivz.
    \end{equation*}
    For a proof, see \citep[Theorem 2.3]{CostabelDauge99}. If $\bfu\in \Hndivz$ solves \eqref{eq:continuous_problem}, taking $\bfv\in\Hndivz\subseteq \X$, we obtain
    \begin{equation*}
        (\bff, \bfv) = a(\bfu, \bfv) = 2(\symgrad(\bfu), \symgrad(\bfv)),
    \end{equation*}
    proving the first claim. The second claim can be proved similarly, since the additional assumption on $\Omega$ implies that $\X$ coincides with $\Hndivz$.
\end{proof}
\Cref{thm:WvsW} has three important consequences: first, non-convexity may lead to
  spurious solutions; second, corners do not require any special treatment; and third, on convex domains with flat
  faces the choice $\alpha = 0$ yields the Stokes solution.  We discuss these three points
  in the remarks below.
  \begin{remark}[Loss of $\mathbf{H}^1$ regularity on
    non-convex domains]
\label{rem:non-convex}
    \Cref{thm:WvsW} does not guarantee that a numerical method based on \eqref{eq:continuous_problem} will converge to the physically correct Stokes solution on non-convex domains.  This observation extends to other $\Hcurl$-based formulations (see, e.g., \citep{HuLeeXu21, LaakmanHuFarrell23, BdVDassiDiPietroDroniou22, DiPietroDroniuQian24, TonnonHiptmairSemiLagrangianNavierStokes, ShipengRuijie, AlonsoRapetti26}). Because these approaches inherently solve a Maxwell-like problem rather than a classical Stokes problem, the resulting velocity field may fail to exhibit full $\mathbf{H}^1(\Omega)$ regularity near re-entrant corners (for both normal and tangential boundary conditions, see \citep[Sect.~2.3]{Amr+98}). This aspect is not always made explicit, and convexity is explicitly assumed only by \cite{DiPietroDroniuQian24}. In \Cref{sec:numerics}, we present an example demonstrating this convergence to a non-physical solution, and we propose a straightforward stabilization technique to recover the correct Stokes behavior.
\end{remark}
\begin{remark}[Curvature and corners] \label{rem:corners}
At corners and edges on the boundary, the normal vector is undefined, but the Weingarten map (or the curvature in two dimensions) can still be defined in a distributional way as a tensor-valued measure on the edges in three dimensions and as a Dirac measure at the corners in two dimensions (see, e.g., \cite{NeunteufelCurvatureComputation}). Consequently, the equivalence between \eqref{eq:NavierSlipBC} and \eqref{eq:sp_d} cannot be understood in a strong, pointwise sense. However, \Cref{thm:WvsW} shows that this equivalence holds strictly in a weak sense on $\Hndivz$, which coincides with $\X$ when $\Omega$ is convex or has a $\C^{1,1}$ boundary, and no singular Dirac terms are required in \eqref{eq:continuous_problem} for consistency. To the best of the authors' knowledge, the fact that it is sufficient to define the Weingarten map on smooth parts of the boundary was first understood by \citet[Theorem 3.1.1.2]{Grisvard}.
\end{remark}
\begin{remark}[Flat faces] \label{rem:flat faces}
As already observed by \citet[Theorem 2.3]{CostabelDauge99}, if $\Omega$ has only flat faces (e.g., a cube), it holds
\begin{equation*}
    (\Curl \bfu, \Curl \bfv) = 2(\symgrad(\bfu), \symgrad(\bfv)), \qquad \forall \bfu, \bfv \in \Hndivz.
\end{equation*}
Therefore, on convex domains, the correct Stokes solution is recovered by taking $\alpha = 0$. As a corollary, the methods considered in \citet{TonnonHiptmairSemiLagrangianNavierStokes, BdVDassiDiPietroDroniou22, DiPietroDroniuQian24, ShipengRuijie} correctly impose slip boundary conditions in general only when $\Omega$ has flat faces. We provide an explicit calculation in \Cref{sec:calculated_example_on_a_cube} showing the consistency with $\alpha = 0$, on the unit cube.
\end{remark}
\begin{corollary}
The kernels of \eqref{eq:continuous_problem} and \eqref{eq:StokesWF} coincide when $\Omega$ is a convex curved polyhedron or when it has a $\C^{1,1}$ boundary.
\end{corollary}
In general, the kernels need not be equal. For example, let $\Omega$ be
a polyhedron with nonzero first Betti number, such as a polygonal annulus in two
dimensions or a solid polyhedral torus in three dimensions. Then $\alpha=0$ almost
everywhere and the kernel of \eqref{eq:continuous_problem} is the nonzero harmonic
space $\h^1$, whereas the kernel of \eqref{eq:StokesWF} is trivial.
}

\section{A priori analysis of the discrete problem}
\label{sec:discrete_vf}
Let $\{\Omega_h\}_h$ be a quasi-uniform and uniformly shape regular family of conforming simplicial meshes on $\Omega$, possibly containing curved elements, with $h > 0$ the mesh width. We assume that the meshes cover $\Omega$ exactly, i.e., $\overline{\Omega}=\bigcup_{K\in\Omega_h}\overline{K}$, that each smooth face of $\Gamma$ is a union of boundary faces, and that each element is the image $K=F_K(\widetilde{K})$ of the reference simplex $\widetilde{K}$ under a $\C^{r+1}$-diffeomorphism $F_K$ such that $\lVert D^sF_K\rVert_{L^\infty(\widetilde{K})}\leq Ch^s$ and $\lVert D^sF_K^{-1}\rVert_{L^\infty(K)}\leq Ch^{-s}$ for $1\leq s\leq r+1$, with $C$ independent of $h$ and $r$ the polynomial degree introduced below, and that neighboring elements induce the same parametrization on a shared face, up to the affine identification of the reference faces. Such meshes can be constructed when $\Gamma$ is piecewise of class $\C^{r+1}$, see \citep[Sect.~6]{Bernardi89}, and affine meshes of polyhedral domains are a particular case. Several results used below are stated in the literature for affine meshes; we indicate briefly why they extend to the present meshes. Let $\Gamma_h$ denote the $(d-1)$-dimensional boundary mesh obtained by restricting $\Omega_h$ to $\Gamma$.

Let $\{\Vh\}_{h>0}$ and $\{\Qh\}_{h>0}$ be families of conforming finite-dimensional subspaces of $\Hcurl$ and $H^1(\Omega)$, respectively, obtained by mapping reference finite element spaces to the physical elements through $F_K$, by composition for $\Qh$ and by the covariant Piola transformation $\widetilde{\bfv}\mapsto(DF_K^{-T}\widetilde{\bfv})\circ F_K^{-1}$ for $\Vh$. We assume that the families are dense in their respective spaces and that $\Grad \Qh \subseteq \Vh$. We further define $
\Qh^0 \coloneqq \left\{q_h \in \Qh \,\middle|\,
\int_{\Omega} q_h\,d\mathbf{x}=0\right\}
$. {More specifically, we take $Q_h$ to be the space of continuous Lagrange polynomials of degree $r$ and $\Vh$ to be the space of N\'{e}d\'{e}lec edge elements of type I \citep{Ned80} of degree $r$, which is denoted by $\mathcal{P}_r^-\Lambda^1$ in \citep{AFW06}, both mapped as above, so that $\Grad\Qh\subseteq\Vh$. In particular, in the lowest order case $r = 1$ we recover the Whitney 1-forms \citep{WhitneyBook}. However, the results we present are general and hold also for other choices of spaces, e.g. \citep{Ned86}, and for cubical meshes \citep{ArnoldAwanou14}, as long as the polynomial degrees are chosen in a compatible way.}

As $\nabla\Qh \subset \Vh$, the pair of spaces satisfies the following inf-sup condition:
\begin{equation} \label{eq: infsup b}
    \inf_{q_h\in \Qh^0}\sup_{\bfv_h\in \Vh}\frac{(\Grad q_h, \bfv_h)}{\lVert q_h\rVert_{H^1(\Omega)}\lVert \bfv_h\rVert_{\Hcurl}} \geq \beta >0,
\end{equation}
with $\beta>0$ independent of $h$. We now define the discrete analogue of $\X$ as
\begin{equation*}
    \X_h \coloneqq \{ \bfv_h \in \Vh \mid (\bfv_h, \Grad q_h) = 0, \, \forall q_h \in \Qh \}.
\end{equation*}

{We introduce a discrete Hodge-Helmholtz orthogonal decomposition for $\X_h$:
\begin{equation*}
    \X_h = \h^1_h \oplus^{\perp} \Z_h,
\end{equation*}
where $\h^1_h\coloneqq \{ \bfv_h \in \X_h \mid \Curl \bfv_h = 0\}$ is the space of discrete harmonics and $\Z_h$ is the $L^2$-orthogonal complement of $\h^1_h$ in $\X_h$. Given a generic function $\bfv_h\in \X_h$ we write $\bfv_h = \bfv_h^{\harm} + \bfv_h^{\perp}$ with $\bfv_h^{\harm}\in \h^1_h$ and $\bfv_h^{\perp}\in\Z_h$.}
With a mild abuse of notation, we also let $a, \widetilde{a}$ and $k$ denote the extension of the bilinear forms from \eqref{eq: bilinear forms} to ${\X + \Vh}$.

The discrete problem reads: find $\bfu_h\in \X_h$ satisfying
\begin{align}
    a(\bfu_h, \bfv_h) &= (\mathbf{f}, \bfv_h)\quad
    \forall \bfv_h \in\X_h.
    \label{eq:discrete_problem}
\end{align}

\subsection{Stability}
Note that $\X_h$ is not a subspace of $\X$ and therefore our approximation is
nonconforming. We can thus not apply the standard theory of conforming approximation of
compact perturbations as described in \citep[Thm.~4.2.9]{SauterSchwabBook} and
\citep[Thm.~13.7]{KressBook}. Instead, our strategy is to apply such results on an
intermediate space $\widehat{\X}_h \subset \X$ and show that the gap between $\X_h$ and
$\widehat{\X}_h$ tends to zero sufficiently fast.  We start by introducing {a key
  technical tool, a curl-preserving lifting.

The lifting operator $L_h:\X_h \to  \X$ is defined as follows: for $\bfv_h = \bfv_h^{\perp} + \bfv_h^{\harm}$ with $(\bfv_h^{\perp}, \bfv_h^{\harm}) \in \Z_h \times \h^1_h$, let
\begin{equation}
	L_h(\bfv_h) = L_h^{\perp}(\bfv_h^{\perp}) + L_h^{\harm}(\bfv_h^{\harm})\in \mathbf{Z} + \mathfrak{H}^1,
\end{equation}
with $L_h^{\perp}:\Z_h \to \Z$ the analogue, for the present boundary conditions, of the Hodge mapping in \citep[Section 4]{Hip02}, defined as the solution of
\begin{align*}
	(\Curl L_h^{\perp}(\bfv_h^{\perp}), \Curl \bfw)  = (\Curl \bfv_h^{\perp}, \Curl \bfw)\quad \forall \bfw \in \mathbf{Z},
\end{align*}
and $L_h^{\harm}:\h^1_h \to \h^1$ is the harmonic lifting defined in \citep{AFW06}:
\begin{align*}
    L_h^{\harm}(\bfv_h^{\harm}) &\coloneqq \bfv_h^{\harm} - \Grad \varphi, \quad
    (\Grad \varphi, \Grad \psi) = (\bfv_h^{\harm}, \Grad \psi), \quad \forall \psi\in H^1(\Omega).
\end{align*}
\citet{AFW06} show that $\lVert L_h^{\harm}(\bfv_h^{\harm}) \rVert_{\VectorLtwo} \leq \lVert \bfv_h^{\harm} \rVert_{\VectorLtwo}$ for each discrete harmonic form $\bfv_h^{\harm}$. The properties of $L_h^{\perp}$ and $L_h^{\harm}$ recalled below are proved for affine meshes in \citep[Lemma~4.5]{Hip02} and \citep[Lemma~5.9]{AFW06}; they extend to the present meshes and boundary conditions, since the interpolation estimate \citep[Lemma~4.6]{Hip02} on which they rely carries over to curved elements by pulling back to $\widetilde{K}$, and the smoothed projection used in \citep[Lemma~5.9]{AFW06} can be replaced by the canonical interpolant $\mathcal{I}_h$ (see the proof of \Cref{lemma:lifting}).
We summarize the properties of $L_h$ in the following lemma.

\begin{lemma}
    \label{lemma:lifting}
    The operator $L_h:\X_h \to \X$ is curl-preserving, injective and $h$-uniformly bounded. Moreover, it satisfies the approximation property
    \begin{equation} \label{eq: lifting props}
        \lVert L_h(\bfv_h) - \bfv_h \rVert_{\Hcurl} \leq Ch^s\lVert L_h(\bfv_h ) \rVert_{\Hcurl} \leq Ch^s\lVert \bfv_h\rVert_{\Hcurl}.
    \end{equation}
\end{lemma}
\begin{proof}
The curl-preserving property follows directly from the definition. The injectivity is due to the injectivity of $L_h^{\perp}$ and $L_h^{\harm}$ \citep[Lemma 5.9]{AFW06}. Next, we show $h$-uniform boundedness by using \Cref{thm:PS} and the $h$-uniform boundedness of $L_h^{\perp}$ and $L_h^{\harm}$:
\begin{align*}
	\lVert L_h(\bfv_h) \rVert_{\Hcurl}
    &\le \lVert L_h^{\perp}(\bfv_h^{\perp})\rVert_{\Hcurl} + \lVert L_h^{\harm}(\bfv_h^{\harm})\rVert_{\Hcurl}\\
    &=  \lVert L_h^{\perp}(\bfv_h^{\perp})\rVert_{\Hcurl} + \lVert L_h^{\harm}(\bfv_h^{\harm})\rVert_{\VectorLtwo}\\
	& \leq \sqrt{1 + C_P^2} \lVert \Curl L_h^{\perp}(\bfv_h^{\perp})\rVert_{L^2(\Omega)^{2d-3}} + \lVert \bfv_h^{\harm} \rVert_{\VectorLtwo}\\
	& = \sqrt{1 + C_P^2}\lVert \Curl \bfv_h \rVert_{L^2(\Omega)^{2d-3}} + \lVert \bfv_h^{\harm} \rVert_{\VectorLtwo}\\
	& \leq \sqrt{1 + C_P^2}\lVert \Curl \bfv_h \rVert_{L^2(\Omega)^{2d-3}} + \lVert \bfv_h \rVert_{\VectorLtwo}
           \leq C \lVert \bfv_h \rVert_{\Hcurl}.
\end{align*}

To show the approximation property, we first recall the approximation property of $L_h^{\perp}$ \citep[Lemma 4.5]{Hip02}:
\begin{align*}
	\lVert L_h^{\perp}\bfv_h^{\perp} - \bfv_h^{\perp} \rVert_{\Hcurl} \leq Ch^s \lVert \Curl \bfv_h^{\perp}\rVert_{L^2(\Omega)^{2d-3}}
	= Ch^s \lVert \Curl L_h^{\perp} \bfv_h^{\perp}\rVert_{L^2(\Omega)^{2d-3}}.
\end{align*}

The analogous result is obtained for $L_h^{\harm}$ by combining \citep[Lemma 5.9]{AFW06} and \Cref{thm:harm} (in \citep[Lemma 5.9]{AFW06}, the smoothed projection of \citep[Theorem 5.6]{AFW06} can be replaced by $\mathcal{I}_h$, because the potential of $\bfv_h^{\harm}-L_h^{\harm}\bfv_h^{\harm}$ is continuous and elementwise in $H^{1+s}$, so that the commuting diagram property applies to it):
\begin{align*}
	\lVert L_h^{\harm}\bfv_h^{\harm} - \bfv_h^{\harm} \rVert_{\Hcurl} &= \lVert
        L_h^{\harm}\bfv_h^{\harm} - \bfv_h^{\harm} \rVert_{\VectorLtwo} \\
	 & \leq Ch^s \lVert L_h^{\harm}\bfv_h^{\harm} \rVert_{\mathbf{H}^s(\Omega)}
	 \leq Ch^s \lVert L_h^{\harm}\bfv_h^{\harm} \rVert_{\VectorLtwo}
	 \leq Ch^s \lVert \bfv_h^{\harm} \rVert_{\VectorLtwo}.
\end{align*}
Putting together these results, we obtain the desired approximation property of $L_h$:
\begin{equation*}
	\begin{split}
	\lVert L_h\bfv_h - \bfv_h \rVert_{\Hcurl} & \leq \lVert L_h^{\perp}\bfv_h^{\perp} - \bfv_h^{\perp}\rVert_{\Hcurl} + \lVert L_h^{\harm}\bfv_h^{\harm} - \bfv_h^{\harm} \rVert_{\Hcurl} \\
	& \leq Ch^s (\lVert \Curl L_h^{\perp} \bfv_h^{\perp} \rVert_{L^2(\Omega)^{2d-3}} + \lVert L_h^{\harm} \bfv_h^{\harm} \rVert_{\VectorLtwo}) \\
	& \leq Ch^s (\lVert \Curl L_h \bfv_h\rVert_{L^2(\Omega)^{2d-3}} + \lVert L_h\bfv_h \rVert_{\VectorLtwo} + \lVert L_h^{\perp} \bfv_h^{\perp} \rVert_{\VectorLtwo})\\
	& \leq Ch^s (\lVert \Curl L_h \bfv_h \rVert_{L^2(\Omega)^{2d-3}} + \lVert L_h \bfv_h \rVert_{\VectorLtwo})\\
	&\leq Ch^s \lVert  L_h\bfv_h \rVert_{\Hcurl}.
	\end{split}
\end{equation*}

Finally, the second inequality in \eqref{eq: lifting props} follows from the $h$-uniform boundedness of $\lift$.
\end{proof}

The injectivity from \Cref{lemma:lifting} implies that $L_h$ admits an inverse on its range $\widehat{\X}_h \coloneqq L_h \X_h$. We continue by showing the properties of this inverse operator $\Pi_h: \widehat{\X}_h \to \X_h$ in the following lemma.
\begin{lemma}
    \label{lemma:elproj}
    The operator $\Pi_h:\widehat{\X}_h \to \X_h$ is curl-preserving, $h$-uniformly bounded and satisfies the following approximation property:
    \begin{equation} \label{eq: inverse lift props}
        \lVert \Pi_h \widehat{\bfv}_h - \widehat{\bfv}_h \rVert_{\Hcurl}
        \leq Ch^s \lVert  \Pi_h\widehat{\bfv}_h\rVert_{\Hcurl}
        \leq Ch^s \lVert \widehat{\bfv}_h\rVert_{\Hcurl}.
    \end{equation}
\end{lemma}
\begin{proof}
The curl-preserving property is straightforward:
\begin{equation*}
	\Curl \Pi_h \widehat{\bfv}_h = \Curl L_h\Pi_h\widehat{\bfv}_h = \Curl \widehat{\bfv}_h.
\end{equation*}
The approximation property of $L_h$ from \Cref{lemma:lifting} implies the $h$-uniform boundedness of $\Pi_h$:
\begin{align*}
	\lVert \Pi_h \widehat{\bfv}_h \rVert_{\Hcurl} & \leq \lVert \Pi_h\widehat{\bfv}_h - L_h\Pi_h\widehat{\bfv}_h \rVert_{\Hcurl} + \lVert  \widehat{\bfv}_h \rVert_{\Hcurl}\\
	&\leq Ch^s\lVert L_h\Pi_h \widehat{\bfv}_h \rVert_{\Hcurl} + \lVert \widehat{\bfv}_h \rVert_{\Hcurl} \\
	& \leq (C \operatorname{diam}(\Omega )^s +1) \lVert \widehat{\bfv}_h \rVert_{\Hcurl}\quad \forall\widehat{\bfv}_h\in \widehat{\X}_h.
\end{align*}
Finally, the approximation property \eqref{eq: inverse lift props} follows from the approximation property of $L_h$ and the stability of $\Pi_h$ shown above:
\begin{equation*}
\begin{split}
	\lVert \Pi_h\widehat{\bfv}_h - \widehat{\bfv}_h \rVert_{\Hcurl} &= \lVert \Pi_h\widehat{\bfv}_h - L_h\Pi_h \widehat{\bfv}_h\rVert_{\Hcurl} \\
    &\leq Ch^s \lVert \Pi_h\widehat{\bfv}_h\rVert_{\Hcurl}
    \leq Ch^s \lVert \widehat{\bfv}_h\rVert_{\Hcurl}\quad\forall\widehat{\bfv}_h\in \widehat{\X}_h.
    \end{split}
\end{equation*}
\end{proof}
}

We proceed by deriving approximation estimates for the operators $\elproj$ and $\lift$ on the boundary $\Gamma$.
\begin{lemma} \label{lem: approx lifting}
There exists a constant $C>0$ independent of $h$ such that the following estimates hold:
\begin{align}
    \lVert \ttrace(\widehat{\bfv}_h - \elproj(\widehat{\bfv}_h))\rVert_{\VectorLtwoGamma} & \leq C h^{s- \frac{1}{2}}{\lVert  \widehat{\bfv}_h\rVert_{\Hcurl}} \quad
    \forall \widehat{\bfv}_h \in \widehat{\X}_h,\label{eq:elproj_bnd} \\
    \lVert \ttrace(\bfv_h - \lift(\bfv_h))\rVert_{\VectorLtwoGamma} &\leq Ch^{s - \frac{1}{2}}{\lVert \bfv_h\rVert_{\Hcurl}} \quad    \forall {\bfv}_h \in {\X}_h\label{eq:lift_bnd}.
\end{align}
\end{lemma}
\begin{proof}
Starting with \eqref{eq:elproj_bnd}, let $\widehat{\bfv}_h\in\widehat{\X}_h$. For each element $K$, put $B_K=DF_K$ and define the local mapped projector by
\begin{align*}
 \widetilde{\bfv}&=B_K^T(\bfv\circ F_K),\qquad
 \overline{\widetilde{\bfv}}
   =|\widetilde{K}|^{-1}\int_{\widetilde{K}}\widetilde{\bfv}\, \mathrm{d\widetilde{x}},\\
 (\boldsymbol{\pi}_{h,0}\bfv)|_K
   &=\bigl(B_K^{-T}\overline{\widetilde{\bfv}}\bigr)\circ F_K^{-1}.
\end{align*}
Its range consists of fields with constant covariant pull-back, so $(\boldsymbol{\pi}_{h,0}\bfv)|_K\in\Vh|_K$. On affine elements it is the $L^2$ projection onto constants; global conformity of $\boldsymbol{\pi}_{h,0}\bfv$ is not required. Let $F$ be a face on $\Gamma$ and let $K_F$ be the unique element to which $F$ belongs. Then we use the triangle inequality to show
\begin{align*}
    \lVert \ttrace( \widehat{\bfv}_h - \elproj ( \widehat{\bfv}_h)) \rVert_{\mathbf{L}^2(F)} & \leq \lVert \widehat{\bfv}_h - \elproj (\widehat{\bfv}_h) \rVert_{\mathbf{L}^2(F)}\\
    & \leq \lVert \widehat{\bfv}_h - \boldsymbol{\pi}_{h,0} \widehat{\bfv}_h\rVert_{\mathbf{L}^2(F)} + \lVert \boldsymbol{\pi}_{h, 0} \widehat{\bfv}_h - \elproj(\widehat{\bfv}_h) \rVert_{\mathbf{L}^2(F)}\\
    & \eqqcolon \mathrm{I} + \mathrm{II}.
\end{align*}
For $\mathrm{I}$, we use the fractional trace inequality \citep[Eq. (12.17)]{EG1} followed by the approximation property of $\boldsymbol{\pi}_{h,0}$ (both are stated for affine meshes and extend to the present meshes by pulling back to $\widetilde{K}$, up to a lower-order term that leads to the full $\mathbf{H}^s$ norm below):
\begin{align*}
    \mathrm{I} &\leq C \left( h_{K_F}^{-\frac{1}{2}} \lVert \widehat{\bfv}_h - \boldsymbol{\pi}_{h,0}\widehat{\bfv}_h \rVert_{\mathbf{L}^2(K_F)} + h_{K_F}^{s - \frac{1}{2}} \lvert \widehat{\bfv}_h - \boldsymbol{\pi}_{h,0}\widehat{\bfv}_h \rvert_{\mathbf{H}^s(K_F)}\right) \\
    & \leq C h_{K_F}^{s - \frac{1}{2}} \lVert \widehat{\bfv}_h \rVert_{\mathbf{H}^s(K_F)}.
\end{align*}
For $\mathrm{II}$, we use the inverse trace inequality \citep[Lem.~12.8]{EG1} (which extends to the present meshes by scaling, since $\boldsymbol{\pi}_{h,0}\widehat{\bfv}_h-\elproj(\widehat{\bfv}_h)\in\Vh|_{K_F}$), followed by a triangle inequality and the approximation properties of $\boldsymbol{\pi}_{h,0}$:
\begin{align*}
    \mathrm{II} & \leq Ch_{K_F}^{- \frac{1}{2}} \lVert \boldsymbol{\pi}_{h,0}\widehat{\bfv}_h - \elproj( \widehat{\bfv}_h )\rVert_{\mathbf{L}^2(K_F)} \\
    & \leq Ch_{K_F}^{-\frac{1}{2}} \left( \lVert \boldsymbol{\pi}_{h,0} \widehat{\bfv}_h - \widehat{\bfv}_h \rVert_{\mathbf{L}^2(K_F)} + \lVert \widehat{\bfv}_h - \elproj(\widehat{\bfv}_h) \rVert_{\mathbf{L}^2(K_F)}\right)\\
    & \leq C\left( h_{K_F}^{s-\frac{1}{2}} \lVert \widehat{\bfv}_h \rVert_{\mathbf{H}^s(K_F)} + h_{K_F}^{-\frac{1}{2}} \lVert \widehat{\bfv}_h - \elproj(\widehat{\bfv}_h) \rVert_{\mathbf{L}^2(K_F)} \right).
\end{align*}
Putting together the bounds for $\mathrm{I}$ and $\mathrm{II}$, and using the quasi-uniformity and shape regularity of the mesh, we obtain the global bound:
\begin{align*}
    \lVert \ttrace(\widehat{\bfv}_h - \elproj ( \widehat{\bfv}_h)) \rVert_{\VectorLtwoGamma} \leq C\left (h^{s- \frac{1}{2}} \lVert \widehat{\bfv}_h\rVert_{\VectorHs} + h^{-\frac{1}{2}}\lVert \widehat{\bfv}_h - \elproj(\widehat{\bfv}_h) \rVert_{\VectorLtwo} \right).
\end{align*}
Then, \eqref{eq:elproj_bnd} follows from \Cref{lemma:elproj} and the embedding inequality $\lVert \widehat{\bfv}_h\rVert_{\VectorHs} \leq C \lVert \widehat{\bfv}_h \rVert_{\Hcurl}$.
To establish \eqref{eq:lift_bnd}, we combine the bound \eqref{eq:elproj_bnd} with the stability of the lifting operator $\lift$:
\begin{align*}
    \lVert \ttrace(\bfv_h - \lift(\bfv_h)) \rVert_{\VectorLtwoGamma}
    &= \lVert \ttrace(\elproj(\lift(\bfv_h)) - \lift(\bfv_h)) \rVert_{\VectorLtwoGamma} \\
    &\leq Ch^{s-\frac{1}{2}}\lVert \lift(\bfv_h)\rVert_{\Hcurl} \\
    &\leq Ch^{s-\frac{1}{2}}\lVert \bfv_h \rVert_{\Hcurl}.
\end{align*}
\end{proof}

As a consequence of \Cref{lem: approx lifting}, we obtain a discrete trace inequality that does not involve any negative power of $h$.
\begin{corollary}
\label{cor:disc_trace_ineq}
There exists a constant $C>0$ independent of $h$ such that the following discrete trace inequality is satisfied for each $h$:
\begin{equation*}
    \lVert \ttrace(\bfv_h)\rVert_{\VectorLtwoGamma}\leq C{\lVert \bfv_h\rVert_{\Hcurl}}, \qquad \forall \bfv_h\in \X_h.
\end{equation*}
\begin{proof}
    The claim follows from a triangle inequality, \eqref{eq:lift_bnd} in \Cref{lem: approx lifting}, \Cref{thm:TraceBounded}, and \Cref{thm:embeddingXinHs}:
    \begin{equation*}
    \begin{split}
      \lVert \ttrace(\bfv_h)\rVert_{\VectorLtwoGamma}&\leq\lVert \ttrace(\bfv_h - \lift(\bfv_h))\rVert_{\VectorLtwoGamma}  + \lVert \ttrace(\lift(\bfv_h))\rVert_{\VectorLtwoGamma}\\
      &\leq C\big(h^{s-\frac{1}{2}}{\lVert \bfv_h\rVert_{\Hcurl}} + {\lVert \lift(\bfv_h)\rVert_{\Hcurl}}\big)\\
      &\leq C\big(h^{s-\frac{1}{2}}{\lVert \bfv_h\rVert_{\Hcurl}} + {\lVert \bfv_h\rVert_{\Hcurl}}\big)\\
      &\leq C \big(\operatorname{diam}(\Omega)^{s-\frac{1}{2}} + 1\big){\lVert
      \bfv_h\rVert_{\Hcurl}}.
    \end{split}
    \end{equation*}
\end{proof}
\end{corollary}

Before presenting the main result, we first analyze the bilinear form $a$ on $\widehat{\X}_h \times \widehat{\X}_h$. For that, we require the following density result.
\begin{lemma} \label{lem:DensityXhHat}
    The family $\{\widehat{\X}_h\}_{h>0}\coloneqq\{\lift(\X_h)\}_{h>0}$ approximates $\X$ in the following limiting sense:
    \begin{equation*}
      \lim_{h\to0}\inf_{\widehat{\bfv}_h\in\widehat{\X}_h}
       \lVert\bfv-\widehat{\bfv}_h\rVert_{\Hcurl}=0
       \qquad\forall\bfv\in\X.
    \end{equation*}
\end{lemma}
\begin{proof}{
   Fix $\epsilon>0$ and let $\bfv\in \X$. We aim to construct a $\bfv_h \in \X_h$ such that $L_h \bfv_h$ is sufficiently close to $\bfv$. First, since $\bC^{\infty}(\overline{\Omega})$ is dense in $\mathbf{H}(\operatorname{curl}, \Omega)$, there exists $\widetilde{\bfv}\in \bC^{\infty}(\overline{\Omega})$ such that
\begin{equation} \label{eq: tilde is epsilon}
	\lVert \bfv - \widetilde{\bfv} \rVert_{\Hcurl} < \epsilon.
\end{equation}
Let $\mathcal{I}_h(\widetilde{\bfv})\in\Vh$ be the standard commuting interpolant of $\widetilde{\bfv}$ satisfying the approximation property (on the present meshes, $\mathcal{I}_h$ is defined by pull-back to $\widetilde{K}$, and the estimate below follows as in \citep[Prop.~3.6]{AylwinJerezHanckes23})
\begin{equation} \label{eq: approx interpolant}
	\lVert \widetilde{\bfv} - \mathcal{I}_h(\widetilde{\bfv})\rVert_{\Hcurl} \leq Ch^{r}( \lVert \widetilde{\bfv} \rVert_{\mathbf{H}^r(\Omega)} + \lVert \Curl \widetilde{\bfv} \rVert_{H^r(\Omega)^{2d-3}}), \qquad r\geq 1.
\end{equation}
Now define $\varphi_h\in Q_h$ via
\begin{equation*}
	(\Grad \varphi_h,\Grad \psi_h) = (\mathcal{I}_h(\widetilde{\bfv}), \Grad \psi_h), \quad \forall \psi_h \in Q_h.
\end{equation*}

We define $\bfv_h \coloneqq \mathcal{I}_h(\widetilde{\bfv}) - \Grad \varphi_h$ and note that it is a viable candidate because it is in $\X_h$. We continue as follows:
\begin{align*}
    \lVert \bfv - L_h \bfv_h \rVert_{\Hcurl}
    &\leq \lVert \bfv - \widetilde{\bfv} \rVert_{\Hcurl} + \lVert \bfv_h - \widetilde{\bfv} \rVert_{\Hcurl} + \lVert (I-L_h) \bfv_h \rVert_{\Hcurl}.
\end{align*}

The first term is bounded by \eqref{eq: tilde is epsilon}. For the second term, we note that
\begin{equation*}
	(\Grad \varphi_h,\Grad  \psi_h) = (\mathcal{I}_h(\widetilde{\bfv}) - \bfv,\Grad  \psi_h),\quad \forall \psi_h \in Q_h,
\end{equation*}
because $\bfv \in \X$. By taking $\psi_h = \varphi_h$, we obtain $\lVert \Grad \varphi_{h} \rVert_{\VectorLtwo}
    \leq \lVert \mathcal{I}_h(\widetilde{\bfv}) - \bfv \rVert_{\VectorLtwo}$. Together with \eqref{eq: approx interpolant}, this leads us to
\begin{equation} \label{eq: second term bound}
\begin{split}
	\| \bfv_h - \widetilde{\bfv} \|_{\Hcurl}
    &\leq \lVert \mathcal{I}_h(\widetilde{\bfv}) - \widetilde{\bfv} \rVert_{\Hcurl} + \| \Grad \varphi_h \|_{\VectorLtwo}\\
    &\leq 2 \lVert \mathcal{I}_h(\widetilde{\bfv}) - \widetilde{\bfv} \rVert_{\Hcurl} + \lVert \widetilde{\bfv} - \bfv\rVert_{\VectorLtwo}\\
    &\leq 2 Ch^r(\lVert \widetilde{\bfv} \rVert_{\mathbf{H}^r(\Omega)} + \lVert \Curl \widetilde{\bfv} \rVert_{H^r(\Omega)^{2d-3}}) + \epsilon.
\end{split}
\end{equation}

For the third term we use the approximation property of $L_h$ from \Cref{lemma:lifting} and \eqref{eq: second term bound}:
\begin{align*}
    \| (I - L_h) \bfv_h \|_{\Hcurl}
    &\le C' h^s \| \bfv_h \|_{\Hcurl} \\
    &\le C' h^s (\| \bfv_h - \widetilde{\bfv} \|_{\Hcurl} + \| \widetilde{\bfv} \|_{\Hcurl}) \\
    &\le C' h^s (2 Ch^r(\lVert \widetilde{\bfv} \rVert_{\mathbf{H}^r(\Omega)} + \lVert \Curl \widetilde{\bfv} \rVert_{H^r(\Omega)^{2d-3}}) + \epsilon
    + \| \widetilde{\bfv} \|_{\Hcurl}) .
\end{align*}

Combining the bounds on the three terms, we obtain
\begin{align*}
    \lVert \bfv - L_h \bfv_h \rVert_{\Hcurl}
    \leq
    \begin{aligned}[t]
      & (2 + C'h^s) \epsilon
    + (1 + C'h^s) 2 Ch^r (\lVert \widetilde{\bfv} \rVert_{\mathbf{H}^r(\Omega)} + \\ & \lVert \Curl \widetilde{\bfv} \rVert_{H^r(\Omega)^{2d-3}})
    + C' h^s \| \widetilde{\bfv} \|_{\Hcurl} .
  \end{aligned}
\end{align*}
Choosing $h$ small enough, we can achieve
\begin{equation*}
	\lVert \bfv - L_h \bfv_h \rVert_{\Hcurl} < 3\epsilon,
\end{equation*}
which concludes the proof.  }
\end{proof}

The results obtained so far allow us to derive the following inf-sup condition on $a$.
\begin{lemma} \label{lem: infsup a on Xhat}
    If {$\ker (a)$} from \eqref{eq: kernel_of_a} is trivial, then there exists $\widehat{h}_0>0$ such that, for each $h\leq \widehat{h}_0$, the bilinear form $a$ satisfies
    \begin{align*}
        \inf_{\widehat{\bfv}_h\in\widehat{\X}_h}
        \sup_{\widehat{\bfw}_h\in\widehat{\X}_h}
        \frac{a(\widehat{\bfv}_h,\widehat{\bfw}_h)}
        {\lVert \widehat{\bfv}_h \rVert_{{\Hcurl}}\lVert \widehat{\bfw}_h \rVert_{{\Hcurl}}} \ge \widehat{\kappa} >0,
    \end{align*}
    with $\widehat{\kappa}$ independent of $h$.
\end{lemma}
\begin{proof}
    We collect four results. First, \Cref{lem:DensityXhHat} shows that $\{\widehat{\X}_h\}_h$ is dense in $\X$. Second, $K:\X\to \X^*$ is compact as shown in the proof of \Cref{thm:FiniteDimensionalKernel}. Third, $\widetilde{a}:\X\times\X\to\R$ is coercive by \Cref{lem:boundednessAndCoercivityBilinearForms}. Finally, $a:\X\times\X\to\R$ has a trivial kernel by assumption. The conditions are now met to invoke \citep[Thm.~4.2.9]{SauterSchwabBook}, from which we obtain the result.
\end{proof}

We are now ready to state and prove the main result of this section.
\begin{theorem} \label{thm:WellPosednessDiscreteProblem}
Assume that {$\ker (a)$}, the null space from \eqref{eq: kernel_of_a}, is trivial. Then, there exists $h_0>0$ and $\kappa>0$ independent of $h$ such that, for each $h\leq h_0$, we have
\begin{equation*}
    \inf_{\bfv_h\in\X_h} \sup_{\bfw_h\in\X_h} \frac{a(\bfv_h,\bfw_h)}{\lVert \bfv_h \rVert_{{\Hcurl}} \lVert \bfw_h \rVert_{{\Hcurl}} } \ge \kappa>0.
\end{equation*}
\end{theorem}
\begin{proof}
{We use a technique similar to that in \cite[Section 7]{HiptmairSchwab02}.}
Given $\bfv\in \X$, let $\widehat{\mathbf{z}}_h(\bfv)\in \widehat{\X}_h$ be the solution of the following problem:
\begin{align} \label{eq: problem zhat}
    a(\widehat{\mathbf{z}}_h(\bfv), \widehat{\mathbf{w}}_h) &= k(\bfv, \widehat{\mathbf{w}}_h), & \forall \widehat{\mathbf{w}}_h &\in \widehat{\X}_h.
\end{align}
Problem \eqref{eq: problem zhat} is well-posed by the continuity and symmetry of $a$ shown in \Cref{lem:boundednessAndCoercivityBilinearForms} and the inf-sup condition from \Cref{lem: infsup a on Xhat}. It now follows from \Cref{lem: infsup a on Xhat} and the continuity of $k$ from \Cref{lem:boundednessAndCoercivityBilinearForms} that
\begin{equation} \label{eq:uniformTestimate}
\begin{aligned}
    \lVert \widehat{\mathbf{z}}_h(\bfv)\rVert_{{\Hcurl}} &\leq \frac{1}{\widehat{\kappa}}\sup_{\widehat{\mathbf{w}}_h\in \widehat{\X}_h}\frac{ a(\widehat{\mathbf{z}}_h(\bfv), \widehat{\mathbf{w}}_h)}{\lVert \widehat{\mathbf{w}}_h\rVert_{{\Hcurl}}} \\
    & = \frac{1}{\widehat{\kappa}}\sup_{\widehat{\mathbf{w}}_h\in \widehat{\X}_h}\frac{  k(\bfv, \widehat{\mathbf{w}}_h) }{\lVert \widehat{\mathbf{w}}_h \rVert_{{\Hcurl}}}
    \leq \frac{C_k}{\widehat{\kappa}} \lVert \bfv\rVert_{{\Hcurl}}.
\end{aligned}
\end{equation}

We define the operator $\hTmap: \widehat{\X}_h \to \widehat{\X}_h$ as
\begin{equation*}
    \hTmap(\widehat{\bfv}_h) \coloneqq \widehat{\bfv}_h - \widehat{\mathbf{z}}_h(\widehat{\bfv}_h),
\end{equation*}
which is $h$-uniformly bounded since $\lVert \hTmap(\widehat{\bfv}_h)\rVert_{{\Hcurl}}\leq (1 + C_k \widehat{\kappa}^{-1})\lVert \widehat{\bfv}_h\rVert_{{\Hcurl}}$ by \eqref{eq:uniformTestimate}.

Next, let $\Tmap  \coloneqq \elproj \circ \hTmap \circ \lift: \X_h \to \X_h$, which is a $h$-uniformly bounded operator as it is the composition of $h$-uniformly bounded operators. We now proceed with the estimate by fixing $\bfv_h\in \X_h$ and setting $\bfw_h \coloneqq \Tmap(\bfv_h)$. By boundedness of $\Tmap$, we immediately obtain
\begin{align} \label{eq: bfw is bounded}
    \lVert \bfw_h \rVert_{{\Hcurl}} \leq C \lVert \bfv_h \rVert_{{\Hcurl}}.
\end{align}

We continue as follows:
\begin{align*}
    a(\bfv_h, \Tmap(\bfv_h)) &= a(\lift(\bfv_h), \Tmap(\bfv_h)) + a(  \bfv_h - \lift(\bfv_h), \Tmap(\bfv_h))\\
    &=  a(\lift(\bfv_h), \hTmap( \lift(\bfv_h))) + a(\lift(\bfv_h), \Tmap(\bfv_h) - \hTmap(\lift(\bfv_h)))\\
    &\quad + a(  \bfv_h - \lift(\bfv_h), \Tmap(\bfv_h))\\
    &\eqqcolon \mathrm{I} + \mathrm{II} + \mathrm{III}.
\end{align*}

Using the shorthand notation $\widehat{\bfv}_h \coloneqq \lift(\bfv_h)$, the first term becomes
\begin{subequations} \label{eqs: bounds}
\begin{equation}
    \begin{split}
        \mathrm{I}
         = a(\widehat{\bfv}_h, \hTmap(\widehat{\bfv}_h))
        &= a(\widehat{\bfv}_h, \widehat{\bfv}_h) - a(\widehat{\bfv}_h, \widehat{\mathbf{z}}_h(\widehat{\bfv}_h)) \\
        &= a(\widehat{\bfv}_h, \widehat{\bfv}_h) - k(\widehat{\bfv}_h, \widehat{\bfv}_h)\\
        &=\widetilde{a}(\widehat{\bfv}_h, \widehat{\bfv}_h)
        {\geq \gamma \lVert \widehat{\bfv}_h \rVert^2_{\Hcurl}}
        {\geq C_{\mathrm{I}}\lVert \bfv_h \rVert^2_{\Hcurl}.}
    \end{split}
\end{equation}
{Note that in the last two steps we have used the coercivity of $\widetilde{a}$ and the $h$-uniform boundedness of $\elproj$.
To bound the second term, we use that $\Pi_h$ is curl-preserving (\Cref{lemma:elproj}) and the approximation property \eqref{eq:elproj_bnd}}:
\begin{equation}
    \begin{split}
    \lvert\mathrm{II}\rvert
    &= \lvert  a(\widehat{\bfv}_h, \elproj(\hTmap(\widehat{\bfv}_h)) - \hTmap(\widehat{\bfv}_h))\rvert \\
    &= \lvert a(\widehat{\bfv}_h, (\elproj - I)\hTmap(\widehat{\bfv}_h))\rvert \\
    & = \lvert \langle \alpha \ttrace(\widehat{\bfv}_h), \ttrace((\elproj - I)\hTmap(\widehat{\bfv}_h))\rangle\rvert\\
    &\leq \lVert \alpha\rVert_{L^{\infty}(\Gamma)}\lVert \ttrace(\widehat{\bfv}_h)\rVert_{\VectorLtwoGamma}\lVert\ttrace((I - \elproj)\hTmap(\widehat{\bfv}_h))\rVert_{\VectorLtwoGamma}\\
    & \leq C h^{s- \frac{1}{2}}\lVert \bfv_h\rVert_{\Hcurl} \lVert \hTmap(\widehat{\bfv}_h)\rVert_{\Hcurl}\\
    &\leq C_{\mathrm{II}} h^{s- \frac{1}{2}}\lVert \bfv_h\rVert_{{\Hcurl}}^2.
    \end{split}
\end{equation}

Finally, we bound the third term in a similar manner, using \eqref{eq:lift_bnd} and Corollary \ref{cor:disc_trace_ineq}:
\begin{equation}
\begin{split}
  \lvert   \mathrm{III}\rvert & = \lvert \langle \alpha\ttrace(\bfv_h - \lift(\bfv_h)), \ttrace(\Tmap(\bfv_h))\rangle \rvert \\
    &\leq \lVert  \alpha \rVert_{L^{\infty}(\Gamma)}\lVert \ttrace(\bfv_h - \lift(\bfv_h))\rVert_{\VectorLtwoGamma} \lVert \ttrace(\Tmap(\bfv_h))\rVert_{\VectorLtwoGamma}\\
    & \leq C h^{s- \frac{1}{2}}\lVert  \bfv_h\rVert_{\Hcurl}\lVert \Tmap(\bfv_h)\rVert_{\Hcurl}\\
    & \leq C_{\mathrm{III}} h^{s- \frac{1}{2}}\lVert \bfv_h\rVert_{{\Hcurl}}^2.
    \end{split}
\end{equation}
\end{subequations}

Collecting \eqref{eqs: bounds}, we have:
\begin{align*}
    a(\bfv_h, \bfw_h) \ge \big({C_{\mathrm{I}}} - h^{s- \frac{1}{2}}(C_{\mathrm{II}} + C_{\mathrm{III}})\big) \lVert \bfv_h\rVert_{{\Hcurl}}^2.
\end{align*}
Since $s > \frac12$, it follows that for $h_0>0$ sufficiently small, there exists $C>0$ such that
\begin{equation} \label{eq: lower bound a}
    a(\bfv_h, \bfw_h) \geq C \lVert \bfv_h\rVert^2_{{\Hcurl}}\quad \forall h\leq h_0.
\end{equation}
Combining \eqref{eq: lower bound a} with \eqref{eq: bfw is bounded} concludes the proof.
\end{proof}

\begin{theorem}[Stability] \label{thm: stability}
    If {$\ker (a)$} is trivial and $h$ is sufficiently small, then the discrete variational problem \eqref{eq:discrete_problem} is well-posed: it admits a unique solution $\bfu_h \in \X_h$ that depends continuously on the data:
    \begin{align*}
        \lVert \bfu_h \rVert_{{\Hcurl}} \leq \frac1{\kappa} \lVert \bff \rVert_{\VectorLtwo}.
    \end{align*}
\end{theorem}
\begin{proof}
    It remains to show that the bilinear form $a$ is continuous on $\X_h \times \X_h$. This follows directly from Cauchy-Schwarz, the boundedness of $\alpha$, and the trace inequality from \Cref{cor:disc_trace_ineq}. Symmetry of $a$ on $\X_h \times \X_h$ follows from the same arguments as in \Cref{lem:boundednessAndCoercivityBilinearForms}. Thus, $a$ is a continuous and symmetric bilinear form that satisfies the inf-sup condition from \Cref{thm:WellPosednessDiscreteProblem}. By the Banach-Ne\v{c}as-Babu\v{s}ka theorem, Problem \eqref{eq:discrete_problem} is well-posed and the result follows.
\end{proof}

\subsection{A priori error estimates}
Throughout this subsection and the following pressure analysis, we assume
$\ker(a)=\{\bfzero\}$ and $0<h\leq h_0$, with $h_0$ chosen as in
\Cref{thm:WellPosednessDiscreteProblem}. Thus the continuous and discrete velocity
solutions are uniquely determined and the discrete inf-sup constant is uniform in $h$.
The consistency identities below remain valid for any existing solutions, independently
of this uniqueness assumption. Whenever a full pressure norm is used, both pressures
are normalized to have zero mean.

We proceed by deriving consistency and error estimates. Since we will be comparing the continuous solution in $\X$ with the discrete solution in $\X_h$, we introduce the space $\X_{\#}\coloneqq \X + \Vh$, endowed with the following norm:
\begin{equation*}
    \lVert \bfv \rVert_{\#}^2 \coloneqq \lVert \bfv\rVert_{\Hcurl}^2 + \lVert \ttrace(\bfv)\rVert_{\VectorLtwoGammaPar}^2.
\end{equation*}
It was shown in \Cref{lem:boundednessAndCoercivityBilinearForms} that $a$ is continuous on $\X \times \X$ and, for sufficiently small $h$, on $\X_h \times \X_h$ in \Cref{thm: stability}. Next, we show that continuity also holds on $\X_{\#} \times \X_h$.

\begin{lemma} \label{lem:boundednessainHastagnorm}
    There exists a constant $C>0$ independent of $h$ such that for all $h$:
    \begin{align*}
        |a(\bfu, \bfv_h)| &\leq C\lVert \bfu\rVert_{\#} \, {
        \lVert \bfv_h\rVert_{\Hcurl}}\quad
        \forall (\bfu, \bfv_h) \in \X_{\#} \times \X_h.
    \end{align*}
\end{lemma}
\begin{proof}
    We have for $\bfu\in\X_{\#}$ and $\bfv_h\in\X_h$
    \begin{align*}
        |a(\bfu,\bfv_h)| &\leq  |(\nabla\times\bfu,\nabla\times\bfv_h)|+|\langle\alpha \ttrace(\bfu),\ttrace(\bfv_h)\rangle|\\
        &\leq \lVert \nabla\times\bfu\rVert_{L^2(\Omega)^{2d-3}} \lVert \nabla\times\bfv_h \rVert_{L^2(\Omega)^{2d-3}} + \lVert\alpha\rVert_{L^\infty(\Gamma)} \lVert\ttrace(\bfu)\rVert_{\VectorLtwoGamma} \lVert\ttrace(\bfv_h)\rVert_{\VectorLtwoGamma}\\
        &\leq C \lVert \bfu\rVert_\# {\lVert \bfv_h\rVert_{\Hcurl}},
    \end{align*}
    where we used \Cref{cor:disc_trace_ineq} to bound $\lVert\ttrace(\bfv_h)\rVert_{\VectorLtwoGamma}$.
\end{proof}

We are now in a position to state the following consistency estimate.
\begin{theorem} \label{thm:GalerkinOrthogonality}
    Let $(\bfu,p)\in (\X\cap \VectorHtwo)\times H^1(\Omega)$ and $\bfu_h\in \X_h$ be solutions of \eqref{eq:strong_problem} and \eqref{eq:discrete_problem}, respectively. Then, $\bfu$ is the solution of \eqref{eq:continuous_problem} and the following consistency estimate holds:
    \begin{equation*}
        \sup_{\bfv_h \in \X_h} \frac{a(\bfu - \bfu_h, \bfv_h)}{\lVert \bfv_h\rVert_{{\Hcurl}}} \leq  \inf_{q_h\in\Qh}\lVert p-q_h\rVert_{H^1(\Omega)}.
    \end{equation*}
\end{theorem}
\begin{proof}
    Since $\bfu\in\VectorHtwo$ solves \eqref{eq:strong_problem}, we have that $-\nvec\times\Curl\bfu+\alpha \ttrace(\bfu)=0$ on $\Gamma$. We thus have for any $\bfv\in\X$
    \begin{equation*}
        \begin{aligned}
            a(\bfu,\bfv) &= (\Curl\bfu,\Curl\bfv) + \langle \alpha \ttrace(\bfu),\ttrace(\bfv)\rangle\\
            &= (\Curl\Curl\bfu,\bfv) + \langle -\nvec\times\Curl\bfu+\alpha \ttrace(\bfu),\ttrace(\bfv)\rangle\\
            &= (\Curl\Curl\bfu,\bfv)\\
            &=(\bff - \Grad p,\bfv)
            = (\bff,\bfv),
        \end{aligned}
    \end{equation*}
    where the last step is due to the orthogonality of $\bfv$ to gradients, cf.~\eqref{eq: def X}. $\bfu$ is therefore the solution of \eqref{eq:continuous_problem}.

    Now, we repeat the same procedure with $\bfv_h\in \X_h$ and $q_h\in\Qh$:
    \begin{align*}
        a(\bfu, \bfv_h) &= (\mathbf{f}, \bfv_h) - (\Grad p, \bfv_h)\\
        &= (\mathbf{f}, \bfv_h) - (\Grad (p - q_h), \bfv_h),
    \end{align*}
    where we now used that $\bfv_h$ is orthogonal to discrete gradients.
    Subtracting \eqref{eq:discrete_problem}, we obtain
    \begin{equation*}
        a(\bfu - \bfu_h, \bfv_h)
        \leq \inf_{q_h\in\Qh}\lVert p-q_h\rVert_{H^1(\Omega)}\lVert \bfv_h\rVert_{{\Hcurl}},
    \end{equation*}
    concluding the proof.
\end{proof}
{
Note that in the case of the Stokes problem \eqref{eq:continuousSymGradStokes}, the regularity assumed in \Cref{thm:GalerkinOrthogonality} holds when $\Gamma$ is $\C^{2,1}$, see \cite[Theorem 1.2]{BdV}. For completeness, we also prove a consistency result for less regular solutions.
\begin{theorem}
Let $\bfu \in \X$ and $\bfu_h\in \X_h$ be solutions of \eqref{eq:continuous_problem} and \eqref{eq:discrete_problem}, respectively. Then the following consistency estimate holds:
\begin{equation*}
    \sup_{\bfv_h \in \X_h}\frac{ a(\bfu - \bfu_h, \bfv_h)}{\lVert \bfv_h \rVert_{\Hcurl}} \leq C\left( h^s \lVert \bff \rVert_{\VectorLtwo} + h^{s - \frac{1}{2}}\lVert \ttrace(\bfu) \rVert_{\mathbf{L}^2(\Gamma)}\right).
\end{equation*}
\end{theorem}
\begin{proof}
Let $\bfv_h \in \X_h$. Using the approximation properties of $\lift$ in \Cref{lemma:lifting} and \Cref{lem: approx lifting}, we obtain
\begin{equation*}
\begin{split}
    a(\bfu, \bfv_h) &= a(\bfu, \lift(\bfv_h)) + a(\bfu, \bfv_h - \lift(\bfv_h))\\
    &= (\bff, \bfv_h) + (\bff, \lift(\bfv_h) - \bfv_h) + \langle \alpha \ttrace(\bfu), \ttrace( \bfv_h - \lift(\bfv_h)) \rangle\\
    & \leq (\bff, \bfv_h) + Ch^s \lVert \bff\rVert_{\VectorLtwo} \lVert \bfv_h \rVert_{\Hcurl}  + Ch^{s- \frac{1}{2}}\lVert \ttrace(\bfu)\rVert_{\mathbf{L}^2(\Gamma)} \lVert \bfv_h\rVert_{\Hcurl}.
    \end{split}
\end{equation*}
Subtracting \eqref{eq:discrete_problem} gives the result.
\end{proof}
}
\Cref{lem:boundednessainHastagnorm} and \Cref{thm:GalerkinOrthogonality} allow us to obtain the following error estimate.
\begin{theorem} \label{thm:OptimalityInfimum}
    If $(\bfu, p) \in (\VectorHtwo\cap\X)\times H^1(\Omega)$ solve \eqref{eq:strong_problem} and $\bfu_h \in \X_h$ solves \eqref{eq:discrete_problem}, then the following error estimate holds
    \begin{equation*}
        \lVert \bfu - \bfu_h\rVert_{\#}\leq C \big(\inf_{\bfv_h\in {\X_h}}\lVert \bfu - \bfv_h\rVert_{\#}+\inf_{q_h\in\Qh}\lVert p-q_h\rVert_{H^1(\Omega)}\big).
    \end{equation*}
\end{theorem}
\begin{proof}
    Let $\bfw_h\in\X_h$. We use a triangle inequality, \Cref{cor:disc_trace_ineq} and \Cref{thm:WellPosednessDiscreteProblem} to derive
    \begin{equation} \label{eq: error1}
        \begin{aligned}
            \lVert \bfu - \bfu_h\rVert_{\#} &\leq \lVert \bfu - \bfw_h \rVert_{\#} + \lVert \bfw_h - \bfu_h \rVert_{\#}\\
            &\leq \lVert \bfu - \bfw_h \rVert_{\#} + C\lVert \bfw_h - \bfu_h \rVert_{{\Hcurl}}\\
            &\leq \lVert \bfu - \bfw_h \rVert_{{\#}} + \frac{C}{\kappa} \sup_{\bfv_h\in\X_h} \frac{ a(\bfw_h-\bfu_h,\bfv_h)}{\lVert\bfv_h\rVert_{{\Hcurl}}}.
        \end{aligned}
    \end{equation}
    We continue by bounding the second term using \Cref{lem:boundednessainHastagnorm} and \Cref{thm:GalerkinOrthogonality}:
    \begin{equation} \label{eq: error2}
        \begin{aligned}
            \sup_{\bfv_h\in\X_h} \frac{ a(\bfw_h-\bfu_h,\bfv_h)}{\lVert\bfv_h\rVert_{{\Hcurl}}}&\leq \sup_{\bfv_h\in\X_h} \frac{ a(\bfw_h-\bfu,\bfv_h)}{\lVert\bfv_h\rVert_{{\Hcurl}}}+\sup_{\bfv_h\in\X_h} \frac{ a(\bfu-\bfu_h,\bfv_h)}{\lVert\bfv_h\rVert_{{\Hcurl}}}\\
            &\leq C \big( \lVert \bfw_h-\bfu\rVert_\# +\inf_{q_h\in\Qh}\lVert p-q_h\rVert_{H^1(\Omega)} \big).
        \end{aligned}
    \end{equation}

    Finally, we combine \eqref{eq: error1} and \eqref{eq: error2} and take the infimum over $\bfw_h$ to obtain the result.
\end{proof}

The error estimate in \Cref{thm:OptimalityInfimum} involves the approximation properties of $\X_h$ in the $\lVert \cdot \rVert_{\#}$ norm{, which we now infer from those of the larger space $\Vh$.}
To do so, we introduce the $L^2$ projection operators $P_{\Grad \Qh}$ and $P_{\X_h}$, onto $\Grad \Qh$ and $\X_h$, respectively. For given $\bfv\in \VectorLtwo$, $P_{\Grad \Qh}\bfv \coloneqq \Grad \varphi_h$ with $\varphi_h$ the solution of
\begin{equation}
    (\Grad \varphi_h, \Grad \psi_h) = (\bfv, \Grad \psi_h) , \qquad \forall \psi_h \in \Qh.
    \label{eq:R_h}
\end{equation}
While $\phi_h$ is only defined up to a constant, $\Grad \varphi_h$ is unique.

These projections have the following important properties:
\begin{subequations} \label{props Pgrad}
\begin{align}
    P_{\Grad \Qh}\bfv &= 0, & \forall \bfv &\in \X, \\
    P_{\Grad \Qh} \bfv_h &= \bfv_h - P_{\X_h}\bfv_h, & \forall \bfv_h &\in \Vh, \\
    \lVert P_{\Grad \Qh} \bfv \rVert_{\VectorLtwo}&\leq \lVert \bfv \rVert_{\VectorLtwo}, & \forall \bfv &\in \VectorLtwo. \label{ineq: l2 pgrad}
\end{align}
\end{subequations}
The $L^2$ stability of $P_{\Grad \Qh}$ in \eqref{ineq: l2 pgrad} follows by taking $\psi_h = \varphi_h$ in \eqref{eq:R_h}.
Remarkably, for $d = 2$, $P_{\Grad \Qh}$ is also stable with respect to the $L^q$ norm with $2 \leq q < \infty$, as shown in the following lemma.

\begin{lemma}
    \label{lem:Lp_stability}
    Let $\Omega \subset \mathbb{R}^2$ be contractible. If $\Gamma$ is of class $\C^\infty$, then there exist constants $C>0$ and $h_1>0$, independent of $h$ and $q$, such that, for all $0<h\leq h_1$,
    \begin{align}
        \lVert P_{\Grad Q_h}\mathbf{v}\rVert_{\mathbf{L}^q(\Omega)}
        &\leq C q \lVert \mathbf{v}\rVert_{\mathbf{L}^q(\Omega)}, &
        \forall \bfv &\in \mathbf{L}^q(\Omega)\quad \forall 2\leq q < \infty.
    \end{align}
\end{lemma}
\begin{proof}
This result is essentially \citep[Corollary~3.1]{Duran}. The only delicate point concerns the $q$-independent $W^{1,q}$-stability of the Neumann Ritz projection, which extends to the meshes considered here by following the arguments in \citep[Sect.~3]{Schatz98}. The assumption $\Gamma\in\C^\infty$ is the one under which the results of \citep{Duran} and \citep{Schatz98} used here are stated.
\end{proof}
\begin{theorem}
\label{thm:from_Xh_to_Vh}
Let $\bfu$ be as in \Cref{thm:OptimalityInfimum} and $h$ sufficiently small. There exists $C >0$ independent of $h$ and of $\bfu$ such that the following approximation property holds:
    \begin{equation} \label{eq: claim 1}
        \inf_{\bfv_h\in\X_h}\lVert \bfu - \mathbf{v}_h\rVert_{\#} \leq \inf_{\bfw_h\in \Vh}\left\{\lVert \bfu - \bfw_h\rVert_{\#} + Ch^{-\frac{1}{2}}\lVert \bfu - \bfw_h\rVert_{\VectorLtwo}\right\}.
    \end{equation}
Moreover, if $d = 2$, $\Omega$ is contractible, $\Gamma$ is of class $\C^\infty$, and $\mathbf{u}\in \mathbf{L}^{\infty}(\Omega)$, the following improved approximation property holds:
\begin{equation} \label{eq: claim 2}
        \inf_{\bfv_h\in\X_h}\lVert \mathbf{u} - \mathbf{v}_h\rVert_{\#}
        \leq \inf_{\bfw_h\in \Vh}\left\{
        \lVert \bfu - \bfw_h\rVert_{\#}
        + {C \lvert \log h\rvert  }
        \lVert \bfu - \bfw_h\rVert_{\mathbf{L}^{\infty}(\Omega)}\right\}.
    \end{equation}
\end{theorem}
\begin{proof}
Using the properties of $P_{\Grad \Qh}$ and $P_{\X_h}$ from \eqref{props Pgrad}, we derive
\begin{align*}
    \inf_{\bfv_h \in \X_h}\lVert \bfu - \bfv_h\rVert_{\#}& \leq \inf_{\bfv_h \in \X_h}\inf_{\bfw_h\in \Vh}\left\{\lVert \bfu - \bfw_h\rVert_{\#} + \lVert \bfw_h - \bfv_h\rVert_{\#}\right\} \\
    &\leq \inf_{\bfw_h\in \Vh}\left\{ \lVert \bfu - \bfw_h\rVert_{\#} + \lVert \bfw_h - P_{\X_h}\bfw_h\rVert_{\#}\right\}\\
    &= \inf_{\bfw_h\in \Vh}\left\{ \lVert \bfu - \bfw_h\rVert_{\#} + \lVert P_{\Grad \Qh}(\bfu - \bfw_h)\rVert_{\#}\right\},
\end{align*}
where we used that $\bfu\in\X$.

We continue by bounding the second term. The $\Hcurl$ part of the $\#$ norm can be bounded by using $\Curl (\Grad \Qh) = 0$ and the stability of $P_{\Grad \Qh}$ from \eqref{ineq: l2 pgrad} as
\begin{equation*}
    \lVert P_{\Grad \Qh} (\bfu - \bfw_h)\rVert_{\Hcurl} = \lVert P_{\Grad \Qh} (\bfu - \bfw_h)\rVert_{\VectorLtwo} \leq \lVert \bfu - \bfw_h\rVert_{\VectorLtwo} .
\end{equation*}

We find an upper bound for the boundary term via a standard inverse inequality \citep[see e.g.][Lem.~12.8]{EG1} (which extends to the present meshes by scaling, since $P_{\Grad\Qh}(\bfu-\bfw_h)\in\Grad\Qh\subseteq\Vh$) and \eqref{ineq: l2 pgrad}:
\begin{equation*}
    \lVert \ttrace(P_{\Grad \Qh}(\bfu - \bfw_h))\rVert_{\VectorLtwoGamma} \leq Ch^{-\frac{1}{2}} \lVert P_{\Grad \Qh}(\bfu - \bfw_h) \rVert_{\VectorLtwo} \leq Ch^{-\frac{1}{2}}\lVert \bfu - \bfw_h\rVert_{\VectorLtwo}.
\end{equation*}
This completes the proof for the first claim \eqref{eq: claim 1}. For the second claim \eqref{eq: claim 2}, we use the {continuous embedding of $L^q(\Gamma)$ in $L^2(\Gamma)$ for $q \ge 2$,} the inverse trace inequality in $L^q$ \citep[Lem.~12.8]{EG1} (again valid on the present meshes by scaling, with a constant independent of $q\geq2$) and the $L^q$-stability of $P_{\Grad \Qh}$ from \Cref{lem:Lp_stability} to obtain
{
\begin{align*}
	\lVert \ttrace(P_{\Grad \Qh}(\bfu - \bfw_h))\rVert_{\mathbf{L}^2(\Gamma)}
        & \leq C \lvert \Gamma\rvert^{-\frac{1}{q}}\lVert\ttrace( P_{\Grad \Qh}(\bfu - \bfw_h))\rVert_{\mathbf{L}^q(\Gamma)}\\
		& \leq  Cq \lvert \Gamma \rvert^{-\frac{1}{q}}h^{-\frac{1}{q}}\lVert \bfu - \bfw_h\rVert_{\mathbf{L}^q(\Omega)}.
\end{align*}
} Shrinking $h_0$ so that $h_0\leq\min\{e^{-2},h_1\}$ and taking $q=\lvert\log h\rvert\geq2$, we obtain the estimate
\begin{equation*}
	\lVert \ttrace(P_{\Grad \Qh}(\bfu - \bfw_h))\rVert_{\VectorLtwoGammaPar} \leq {C \lvert \log h\rvert }\lVert \bfu - \bfw_h\rVert_{\mathbf{L}^{\infty}(\Omega)},
\end{equation*}
Here $h^{-1/\lvert\log h\rvert}=e$ and
$\lvert\Gamma\rvert^{-1/\lvert\log h\rvert}
\leq\max\{1,\lvert\Gamma\rvert^{-1/2}\}$ for $0<h\leq e^{-2}$.
The factor $\lvert\Omega\rvert^{1/q}$ from the embedding of
$L^\infty(\Omega)$ into $L^q(\Omega)$ is uniformly bounded as well.
This concludes the proof.
\end{proof}

{
\begin{remark}The improved error estimates rely on the $L^q$-stability of the projection $P_{\Grad Q_h}$. To the best of the authors' knowledge, this result has been established in the literature only for topologically trivial planar domains. We conjecture that the stability property extends to three dimensions as well, as our numerical experiments do not reveal any deterioration in the observed convergence rates; see Section~\ref{sec:numerics}. \end{remark}}

We conclude this section by commenting on the approximation estimates of $\Vh$ in the
$\left\|\cdot\right\|_{\#}$
norm.
{Recall that $\mathcal{I}_h:\mathbf{H}^2(\Omega) \to \Vh$ is the standard interpolation operator already introduced in the proof of \Cref{lem:DensityXhHat}. With a slight abuse of notation, we also denote by $\mathcal{I}_h$ the
interpolator on the boundary trace space, $ \mathcal{I}_h: \ttrace \mathbf{H}^2(\Omega) \to \ttrace \Vh$.
On the present meshes, $\mathcal{I}_h$ is defined by pull-back to $\widetilde{K}$; it still commutes with $\ttrace$, and its approximation properties extend to curved elements by the same pull-back argument, see \citep[Prop.~3.6]{AylwinJerezHanckes23} for the $\mathbf{H}(\operatorname{curl})$-estimate. To estimate boundary terms, use the commutativity of $\ttrace$ with $\mathcal{I}_h$:
\begin{equation*}
    \lVert \ttrace( \bfu - \mathcal{I}_h(\bfu))\rVert_{\mathbf{L}^2(\Gamma)} = \lVert \ttrace(\bfu) - \mathcal{I}_h(\ttrace(\bfu))\rVert_{\mathbf{L}^2(\Gamma)} \leq Ch^r \lVert \ttrace(\bfu) \rVert_{\mathbf{H}^r_{\mathrm{pw}}(\Gamma)}.
\end{equation*}
}

Therefore, taking $\bfw_h = \mathcal{I}_h(\bfu)$ in \Cref{thm:from_Xh_to_Vh} and $q_h$ to be the Lagrange interpolant of $p$ in \Cref{thm:OptimalityInfimum} (see \citep[Lemma~3.2]{LiQiuXieYu24} for exactly fitted curved meshes), we obtain the following error estimates.
\begin{theorem} \label{thm: optimal rates}
{Let $(\bfu,p)$ and $\bfu_h$ be the continuous and discrete
solutions specified in \Cref{thm:OptimalityInfimum}, under the standing
stability and pressure-normalization assumptions of this subsection.
Assume that $\bfu \in \mathbf{H}^2(\Omega) \cap \mathbf{H}^r(\Omega)$,
$\ttrace(\bfu)\in\mathbf{H}^r_{\mathrm{pw}}(\Gamma)$,
$\Curl\bfu\in \mathbf{H}^r(\Omega)$ and $p\in H^{r+1}(\Omega)$.
Then the following error estimate holds:
\begin{equation*}
    \lVert \bfu - \bfu_h \rVert_{\#} \leq C h^{r - \frac{1}{2}}\left( \lVert \bfu
      \rVert_{\mathbf{H}^r(\Omega)} + \lVert \Curl \bfu \rVert_{H^r(\Omega)^{2d-3}} +
      \lVert \ttrace(\bfu ) \rVert_{\mathbf{H}^r_{\mathrm{pw}}(\Gamma)} + \lVert p
      \rVert_{H^{r+1}(\Omega)} \right).
\end{equation*}
Moreover, if $d = 2$, $\Omega$ is contractible, $\Gamma$ is of class $\C^\infty$, and $\bfu\in\mathbf{W}^{r, \infty}(\Omega)$, the following estimate, optimal up to a logarithmic factor, holds:
\begin{equation*}
  \lVert \bfu - \bfu_h \rVert_{\#} \leq C h^{r}\bigl(
  \begin{aligned}[t]
    & \lVert \bfu \rVert_{\mathbf{H}^r(\Omega)} + \lVert \Curl \bfu
    \rVert_{H^r(\Omega)^{2d-3}} + \lVert \ttrace(\bfu ) \rVert_{\mathbf{H}^r_{\mathrm{pw}}(\Gamma)} +
    \\
    & \lvert\log h\rvert \lVert \bfu \rVert_{\mathbf{W}^{r,\infty}(\Omega)} +  \lVert p \rVert_{H^{r+1}(\Omega)} \bigr).
  \end{aligned}
      \end{equation*}}
\end{theorem}
\begin{remark}
    For every $0<\epsilon<r$,
$\lvert\log h\rvert\leq C(\epsilon)h^{-\epsilon}$ for sufficiently small $h$;
see \citep[Lem.~10]{Nitsche}.
Consequently, the improved estimate in \Cref{thm: optimal rates} gives
an error of order $h^{r-\epsilon}$, with a constant depending on $\epsilon$.
This is arbitrarily close to order $r$, which is observed in the smooth
numerical examples of \Cref{sec:numerics}.
\end{remark}

\subsection{Saddle-point formulation and error estimates for the pressure}

A direct implementation of the space $\X_h$ is impractical. Instead, we solve the
following saddle-point problem: find $(\bfu_h, p_h)\in \Vh\times \Qh^0$ satisfying
\begin{subequations}
    \label{eq:discrete_saddlepoint}
    \begin{align}
        a(\bfu_h, \bfv_h) + b(\bfv_h, p_h) &= (\mathbf{f}, \bfv_h) & \forall \bfv_h &\in \Vh, \\
        b(\bfu_h, q_h) &= 0 & \forall q_h &\in \Qh^0,
    \end{align}
\end{subequations}
with $b(\bfv_h, q_h)\coloneqq (\bfv_h, \Grad q_h)$. Let $B: \Vh \to (\Qh^0)^*$ be the associated linear operator.
\begin{lemma}
    Problem \eqref{eq:discrete_saddlepoint} is well-posed if and only if problem \eqref{eq:discrete_problem} is well-posed.
\end{lemma}
\begin{proof}
Following standard saddle-point theory, Problem \eqref{eq:discrete_saddlepoint} is well-posed if $a$ satisfies an inf-sup condition on
\begin{equation*}
    \ker B\coloneqq \{ \bfv_h \in \Vh\mid b(\bfv_h, q_h) = 0,\, \forall q_h\in \Qh^0\},
\end{equation*}
and $b$ satisfies an inf-sup condition. The former condition is equivalent to the well-posedness of \eqref{eq:discrete_problem} since $\X_h = \ker B$, while the latter is satisfied by assumption \eqref{eq: infsup b}.
\end{proof}

The equivalence between the problems is summarized in the following lemma, which we present without proof.
\begin{lemma}
Let $(\bfu_h, p_h)\in \Vh\times \Qh^0$ be a solution of \eqref{eq:discrete_saddlepoint}. Then $\bfu_h$ solves \eqref{eq:discrete_problem}. Vice versa, if $\bfu_h\in \X_h$ solves \eqref{eq:discrete_problem}, then $(\bfu_h, p_h)$ solves \eqref{eq:discrete_saddlepoint} where $p_h\in \Qh^0$ is the unique solution of
\begin{equation*}
    (\Grad p_h, \Grad q_h) = (\mathbf{f}, \Grad q_h) - a(\bfu_h, \Grad q_h) \qquad \forall q_h\in \Qh^0.
\end{equation*}
\end{lemma}

\Cref{cor:disc_trace_ineq} does not provide an $h$-uniform
trace bound for arbitrary test functions in $\Vh$.
To analyze \eqref{eq:discrete_saddlepoint}, we therefore use the following
mesh-dependent estimate.

\begin{lemma} \label{lem: a_bound_sharp_gradients}
    Let $\bfu\in \X_{\#}$ and $\bfv_h \in \nabla Q_h$. Then there exists a constant $C_a > 0$ such that
    \begin{align}
        \lvert a(\bfu, \bfv_h)\rvert \leq C_a h^{-\frac12} \lVert \bfu\rVert_{\#} \lVert \bfv_h \rVert_{\VectorLtwo}.
    \end{align}
\end{lemma}
\begin{proof}
We use $\nabla \times \bfv_h = 0$ with a Cauchy-Schwarz inequality and a standard inverse trace inequality on the discrete $\bfv_h$ (see e.g.~\citep[Lem.~12.8]{EG1}; it extends to the present meshes by scaling):
\begin{gather*}
  \begin{aligned}
    \lvert a(\bfu, \bfv_h)\rvert = \bigl|\langle \alpha \ttrace( \bfu), \ttrace(\bfv_h) \rangle\bigr| & \leq \lVert \alpha\rVert_{L^\infty(\Gamma)} \lVert
    \ttrace(\bfu)\rVert_{\VectorLtwoGammaPar} \lVert
    \ttrace(\bfv_h)\rVert_{\VectorLtwoGammaPar} \\ & \leq C_a \lVert \bfu\rVert_{\#}
    (h^{-\frac12} \lVert \bfv_h \rVert_{\VectorLtwo}).
  \end{aligned}
\end{gather*}
\end{proof}
This yields the following error estimate for the pressure.
\begin{theorem} \label{thm: error estimate p}
Let $(\bfu,p)\in(\X\cap\VectorHtwo)\times
\bigl(H^1(\Omega)\cap L^2_0(\Omega)\bigr)$ solve
\eqref{eq:strong_problem}, and let
$(\bfu_h,p_h)\in\Vh\times Q_h^0$ solve
\eqref{eq:discrete_saddlepoint}, under the standing stability assumptions.
Then the following estimate holds:
    \begin{equation*}
    \lVert p - p_h\rVert_{H^1(\Omega)} \leq C \left(h^{-\frac{1}{2}} \lVert \bfu - \bfu_h\rVert_{\#} + \inf_{q_h\in \Qh}\lVert p - q_h\rVert_{H^1(\Omega)}\right).
\end{equation*}
\end{theorem}
\begin{proof}
Since $(\bfu, p)$ is the strong solution, Galerkin orthogonality holds:
\begin{align}
    a(\bfu - \bfu_h, \bfv_h) + b(\bfv_h, p - q_h) &= b(\bfv_h, p_h - q_h), &
    \forall (\bfv_h, q_h) &\in \Vh \times \Qh^0.
    \label{eq: Galerking orth}
\end{align}

Using the inf-sup condition on $b$ from \eqref{eq: infsup b}, \eqref{eq: Galerking orth}, and \Cref{lem: a_bound_sharp_gradients}, we derive for any $q_h \in Q_h^0$:
\begin{align*}
    \beta \lVert p_h - q_h\rVert_{H^1(\Omega)}
    &\leq \sup_{\bfv_h\in \Vh}\frac{ b(\bfv_h, p_h - q_h) }{\lVert \bfv_h\rVert_{\Hcurl}}
    = \sup_{\bfv_h\in \Grad\Qh}\frac{ b(\bfv_h, p_h - q_h) }{\lVert \bfv_h\rVert_{\VectorLtwo}}\\
    &= \sup_{\bfv_h\in \Grad\Qh}\frac{ a(\bfu - \bfu_h, \bfv_h) + b(\bfv_h, p - q_h) }{\lVert \bfv_h\rVert_{\VectorLtwo}}\\
    &\leq  \sup_{\bfv_h\in \Grad\Qh}\frac{C_a  h^{-\frac{1}{2}} \lVert\bfu - \bfu_h\rVert_{\#}\lVert \bfv_h\rVert_{\VectorLtwo} +  \lVert p - q_h\rVert_{H^1(\Omega)}\lVert  \bfv_h\rVert_{\VectorLtwo}}{\lVert \bfv_h\rVert_{\VectorLtwo}} \\
    &= C_a h^{-\frac{1}{2}}\lVert\bfu - \bfu_h\rVert_{\#} +  \lVert p - q_h\rVert_{H^1(\Omega)}.
\end{align*}

The result now follows by using the triangle inequality $\lVert p - p_h\rVert_{H^1(\Omega)}\leq \lVert p - q_h\rVert_{H^1(\Omega)} + \lVert p_h - q_h\rVert_{H^1(\Omega)}$ and taking the infimum over $q_h\in Q_h^0$.
This infimum equals the one over $Q_h$ in the stated estimate:
for mean-zero $p$, replacing $q_h$ by
$q_h^0=q_h-|\Omega|^{-1}\int_\Omega q_h\, \mathrm{d} \mathbf{x}$ leaves its gradient unchanged and gives
\begin{equation*}
 \lVert p-q_h^0\rVert_{H^1(\Omega)}^2
 =
 \lVert p-q_h\rVert_{H^1(\Omega)}^2
 -
 |\Omega|^{-1}\left|\int_\Omega q_h\, \mathrm{d} \mathbf{x}\right|^2.
\end{equation*}
\end{proof}

Combining this estimate with the velocity bound gives
an upper error bound with a half-order lower power of $h$.

\section{Numerical experiments}
\label{sec:numerics}

In this section, we conduct various numerical experiments to validate our theoretical
results and demonstrate the practical relevance of the method. Specifically, we { examine
  two of the tests performed in \citep{LimacheObjectivityTests}. We compare
  our method against the} two $H^1$-based formulations considered therein, which are
referred to as the {\lq\lq Divergence\rq\rq}{} and the {\lq\lq Laplace\rq\rq}{}
formulations. In the case of straight boundaries, these two methods converge to the same
solution. However, in the case of curved boundaries, the solutions of these methods
differ. Only the \lq\lq Divergence\rq\rq\, formulation satisfies the objectivity
requirement.  In both tests, our method (which will be referred to as
\lq\lq $H(\mathrm{curl})$\rq\rq\, in the figures) produces the solution required by the objectivity criteria, as does the \lq\lq
Divergence\rq\rq\, formulation on curved meshes.

Unless stated otherwise, we use third-order N\'ed\'elec edge elements of the first kind to discretize the velocity and standard $H^1(\Omega)$-conforming elements of the same order to discretize the pressure.
On smooth domains, we approximate the Weingarten map using the method described in \citep[Sect.~3]{NeunteufelCurvatureComputation}. This approach allows us to obtain the correct solution even when a non-curved mesh is used, where traditional methods may fail, see \Cref{sec:Annulus}. However, to reduce the influence of the geometry and Weingarten-map approximations on the convergence study (the analysis assumes exactly fitted curved elements), we employ curved meshes in the convergence tests. For lowest-order (first-order) finite element spaces, the curved boundary is approximated using polynomials of degree three. More generally, for finite element spaces of order \(r\), we use polynomials of degree \(r+2\) to represent the curved elements. We note that on nonsmooth domains the method of \citep[Sect.~3]{NeunteufelCurvatureComputation} cannot be applied, since it produces a Riesz representative of the distributional Weingarten map, which may become unbounded at corners and edges. In fact, our analysis requires \(\alpha\) to be bounded in \(\mathbf{L}^\infty(\Gamma)\), and we have shown in \Cref{thm:WvsW} and \Cref{rem:corners} that, under the assumptions of the theorem, it suffices to take \(\alpha = 2\mathcal{W}\) on the smooth parts of the boundary (\(\alpha = 0\) on flat faces) and to neglect the singular part of the distributional Weingarten map on edges and vertices. In the case of Dirichlet boundary conditions, we apply the Nitsche-type method of \citet{Nitsche}, as described in \Cref{sec: Nitsche}.
All numerical tests are performed with the open-source finite element library
\texttt{NGSolve} \citep{ngSolve}. The source code is freely available at
\url{https://gitlab.com/WouterTonnon/vvphcurlslip}.

\subsection{A manufactured solution in 2D}
\label{sec:experiment1}

\begin{figure}[h]
  \centering
  \includegraphics[width = 0.49\textwidth]{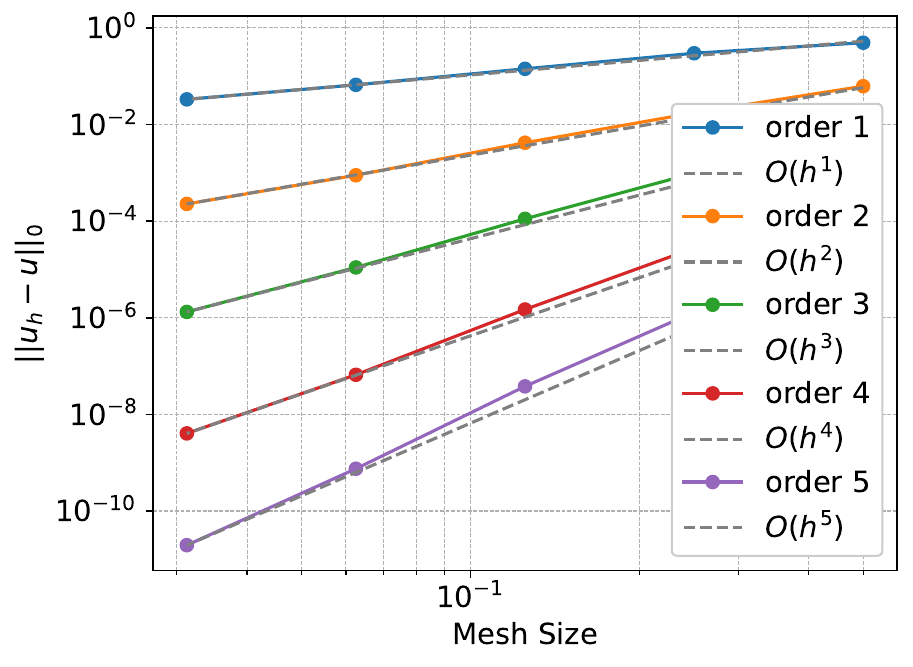}
  \includegraphics[width = 0.49\textwidth]{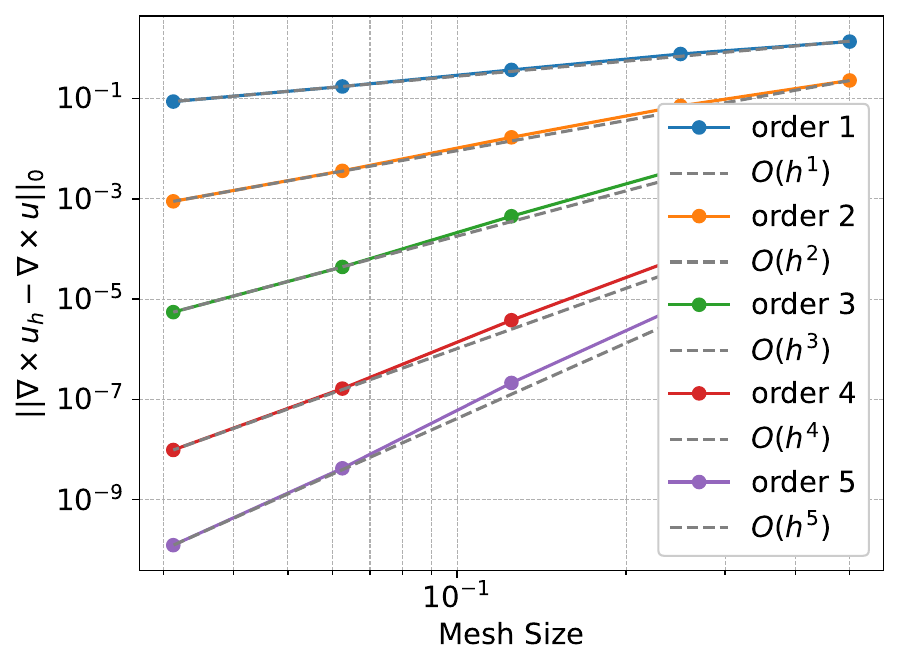} \\
  \includegraphics[width = 0.49\textwidth]{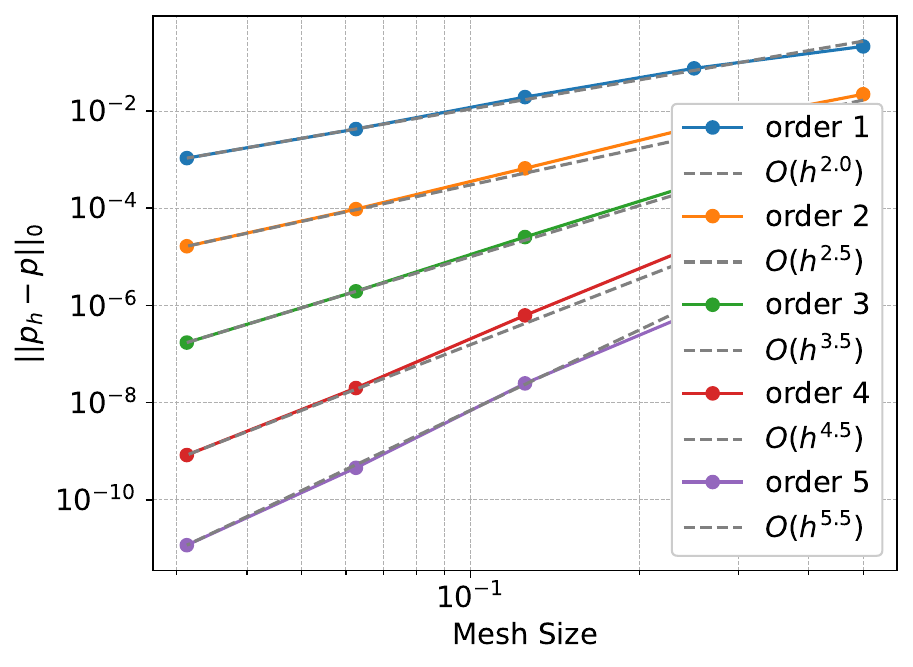}
  \includegraphics[width = 0.49\textwidth]{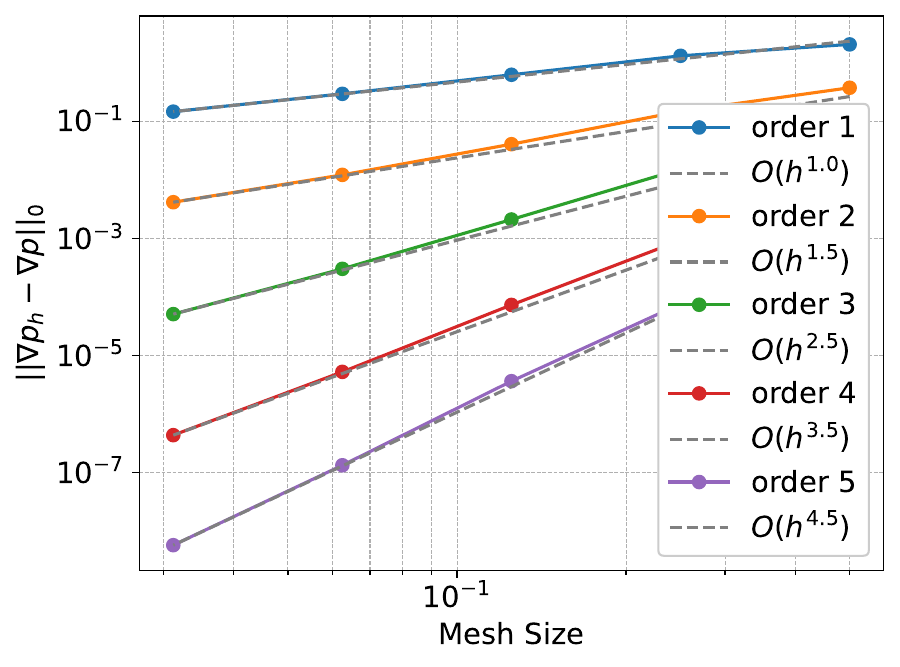}
  \caption{Convergence analysis of (top) $\bfu$ and $\nabla\times\bfu$ in the $L^2$ norm, and (bottom) $p$ and $\nabla p$ in the $L^2$ norm for the experiment as discussed in \cref{sec:experiment1}. The results for lowest-order elements are labeled \lq\lq order 1\rq\rq.
  }
  \label{fig:RobinConvergenceAnalysis}
\end{figure}

In this experiment, we consider a manufactured solution on an ellipse $\Omega\subset \R^2$
with width 4 and height 2 to validate the convergence rates for the velocity and
pressure. We seek $\bfu\in \bC^\infty(\Omega)$ and $p\in \C^\infty(\Omega)$ such that
\begin{equation} \label{eq:alphaEquation}
    \begin{aligned}
        \nabla\times\nabla\times\bfu + \nabla p &= \bff,&\text{ in }\Omega,\\
        \nabla\cdot\bfu &= 0,&\text{ in }\Omega, \\
        \bfu\cdot\nvec &= z,&\text{ on }\Gamma, \\
        -\nvec\times\nabla\times\bfu + (2\mathcal{W})\nvec\times(\bfu\times\nvec) &=\bfg,&\text{ on }\Gamma,
    \end{aligned}
\end{equation}
where $\mathcal{W}\in L^\infty(\Gamma)$ is the Weingarten map (a scalar for
$\Omega\subset \R^2)$. We choose the given functions $\bff$, $z$, and $\bfg$ such that
\begin{align}
    \bfu &= \begin{bmatrix}
        -\sin(2x)\cos(2y)\\
        \cos(2x)\sin(2y)
    \end{bmatrix}, &
    p = x\sin(3x){\cos(y)},
\end{align}
is the solution of the boundary value problem \eqref{eq:alphaEquation}. Even though the
curvature of an ellipse is known exactly, we emphasize that an exact curvature is often
either not available or difficult to represent numerically. For that reason, we will
approximate the curvature as described in the beginning of this section. In
\Cref{fig:RobinConvergenceAnalysis}, we display the $L^2$-error of $\bfu$, the $L^2$-error of $\nabla\times\bfu$, the $L^2$-error of $p$, and the $L^2$-error of $\nabla p$,
respectively. All pressure errors are computed after normalizing both pressures to have zero mean. The convergence rates are in agreement with our analysis from \Cref{thm:
  optimal rates}. Moreover, for $r\geq2$ we observe half an order loss in convergence for the $L^2$ norm of $\nabla p$, consistent with the upper bound of \Cref{thm: error estimate p}, while for $r=1$ the observed rate $1$ is better than this bound.

\subsection{A non-smooth solution} \label{sec:NonSmoothExperiment}
\begin{figure}[h!]
  \centering
  \begin{minipage}[t]{0.49\textwidth}
    \centering
    \includegraphics[width=\linewidth]{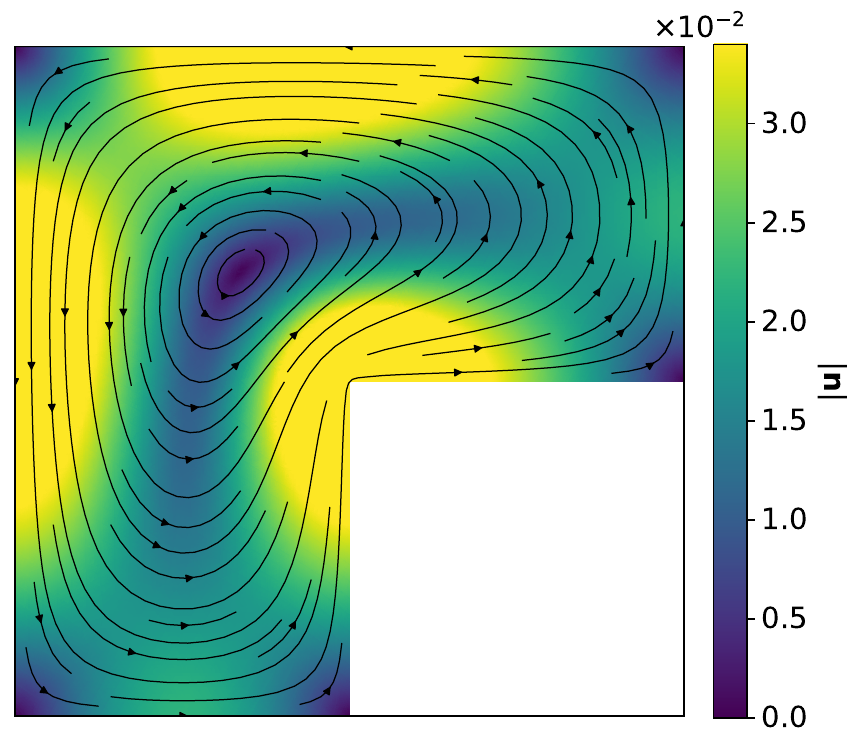}
    \\[-0.5ex]{\small (a) Standard $\Hcurl$ formulation}
  \end{minipage}\hfill
  \begin{minipage}[t]{0.49\textwidth}
    \centering
    \includegraphics[width=\linewidth]{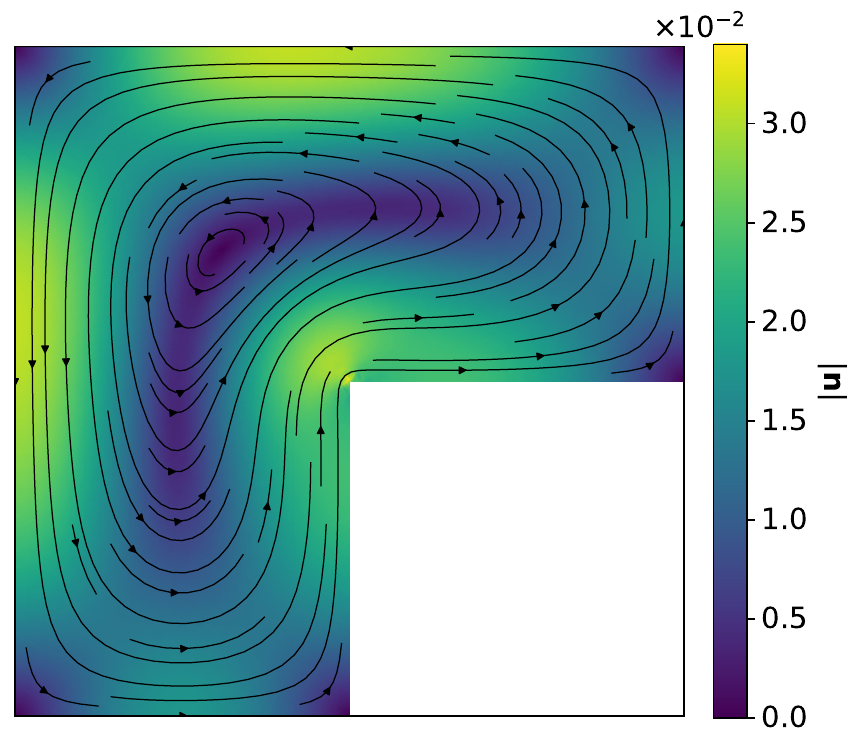}
    \\[-0.5ex]{\small (b) Stabilized $\Hcurl$ formulation}
  \end{minipage}

  \vspace{0.75ex}

  \begin{minipage}[t]{0.49\textwidth}
    \centering
    \includegraphics[width=\linewidth]{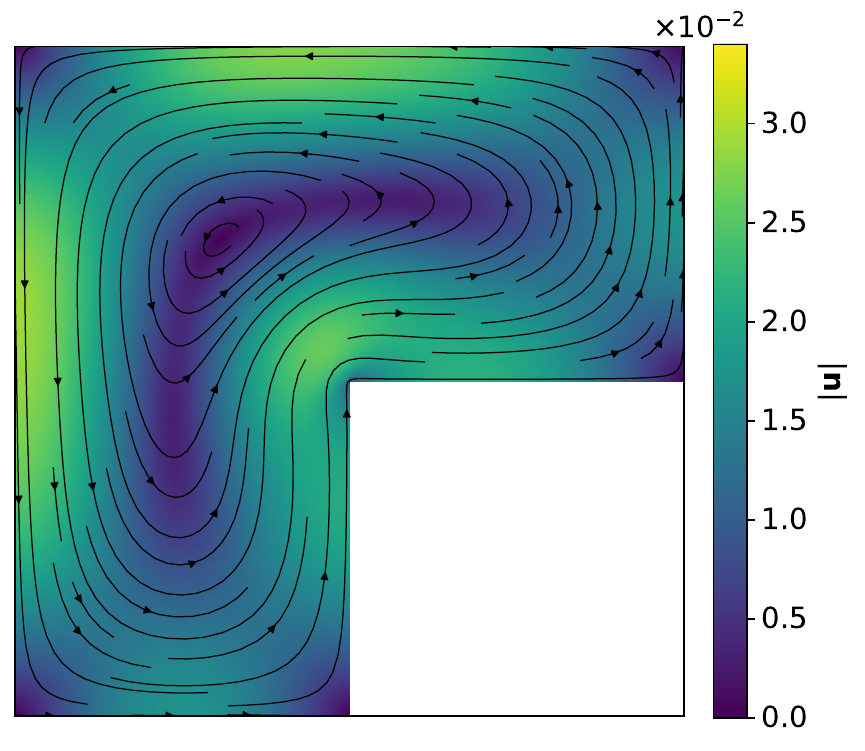}
    \\[-0.5ex]{\small (c) Reference ``Divergence'' solution}
  \end{minipage}\hfill
  \begin{minipage}[t]{0.49\textwidth}
    \centering
    \includegraphics[width=\linewidth]{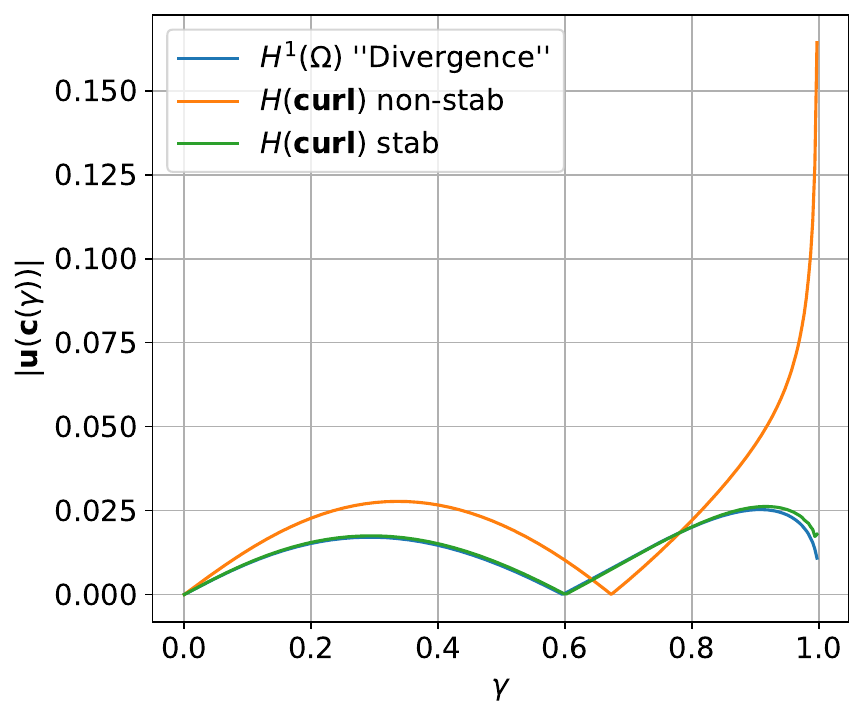}
    \\[-0.5ex]{\small (d) Line comparison of the velocity field}
  \end{minipage}

  \caption{Visual comparison of the computed velocity field of the $\Hcurl$ formulation with and without stabilization, and comparison to the reference ``Divergence'' solution described in \cref{sec:NonSmoothExperiment}. Subplot (d) plots the magnitude of the velocity of the different velocity fields over the curve $\mathbf{c}(\gamma)\coloneqq (1-\gamma)[-1,1]^T$, $\gamma\in [0,1]$. For all computations, we used $h=0.015$ and $C_w=10^4$.}
  \label{fig:LshapedStabilizedCompare}
\end{figure}

First, we consider an experiment for which the exact solution is not known. We take the L-shaped domain $\Omega = (-1,1)^2 \setminus  ([0,1)\times (-1,0] )$, with $\bff = [1,x]^T$, and impose homogeneous Navier slip boundary conditions on the entire boundary.
Since all boundary segments are straight, we set $\alpha = 0$, cf.\ \Cref{rem:flat faces}. As a reference solution we use the solution obtained with the ``Divergence'' formulation, which is known to converge to the exact Stokes solution. In contrast, the $\Hcurl$ formulation fails to converge to the physical solution; see \Cref{fig:LshapedStabilizedCompare}(a). This behavior is consistent with \Cref{thm:WvsW}: on non-convex domains, the weak formulations \eqref{eq:continuous_problem} and \eqref{eq:StokesWF} are in general not equivalent, and consequently the corresponding conforming numerical methods may converge to different limits, cf.\ \Cref{rem:non-convex}.

To recover the correct Stokes solution, we propose a simple remedy. We augment the discrete bilinear form with the stabilization term
\begin{equation}
\frac{C_w}{h}\sum_{f\in\mathring{\mathcal{F}}_h} \int_f \jump{\bfu_h} \cdot \jump{\bfv_h}\,\mathrm{dS},
\label{eq:stab}
\end{equation}
where $\mathring{\mathcal{F}}_h$ denotes the set of all interior facets, $\jump{\cdot}$ is the standard jump operator, and $C_w$ is a user-defined constant.

The resulting numerical solution is shown in \Cref{fig:LshapedStabilizedCompare}(b) and agrees well with the reference solution in \Cref{fig:LshapedStabilizedCompare}(c). This is further illustrated by the line comparison in \Cref{fig:LshapedStabilizedCompare}(d), where we plot the velocity magnitude along the curve $\mathbf{c}(\gamma)=(1-\gamma)[-1,1]^T$, $\gamma\in[0,1]$. In this experiment, the jump penalization \eqref{eq:stab} appears to restore the convergence towards the Stokes solution also on this non-convex domain.

Next, we consider a manufactured solution on the same domain, following \citep{HoustonSchoetzau,verfuhrt1996review}. Let $(r,\phi)$ denote the standard polar coordinates, then we seek the solution
\begin{align}
    \bfu &=r^\lambda \begin{bmatrix}
        (1+\lambda)\sin(\phi)\Psi(\phi)+\cos(\phi)\Psi'(\phi)\\
        \sin(\phi)\Psi'(\phi)-(1+\lambda)\cos(\phi)\Psi(\phi)
    \end{bmatrix}, \\
    p &= -r^{\lambda-1}[(1+\lambda)^2\Psi'(\phi)+\Psi'''(\phi)]/(1-\lambda),
\end{align}
where
\begin{align}
    \Psi(\phi) &= \sin((1+\lambda)\phi)\cos(\lambda\omega)/(1+\lambda)-\cos((1+\lambda)\phi)\notag\\&-\sin((1-\lambda)\phi)\cos(\lambda\omega)/(1-\lambda)+\cos((1-\lambda)\phi),\\
    \omega &=\frac{3\pi}{2},
\end{align}
and $\lambda\approx 0.54448373678246$. Note that both $\nabla\bfu$ and $\nabla p$ are singular at the origin, in particular $\bfu\not\in\VectorHtwo$ and $p\notin H^1(\Omega)$. We take $\bff=\mathbf{0}$ and impose Dirichlet boundary conditions on the part of the boundary that coincides with the boundary of $(-1,1)^2$ and inhomogeneous Navier slip boundary conditions on the remaining part of the boundary, with data given by the exact solution. Since all boundaries are straight, the curvature is set to zero, cf. \Cref{rem:flat faces}.

\begin{figure}[h]
  \centering
  \includegraphics[width = 0.49\textwidth]{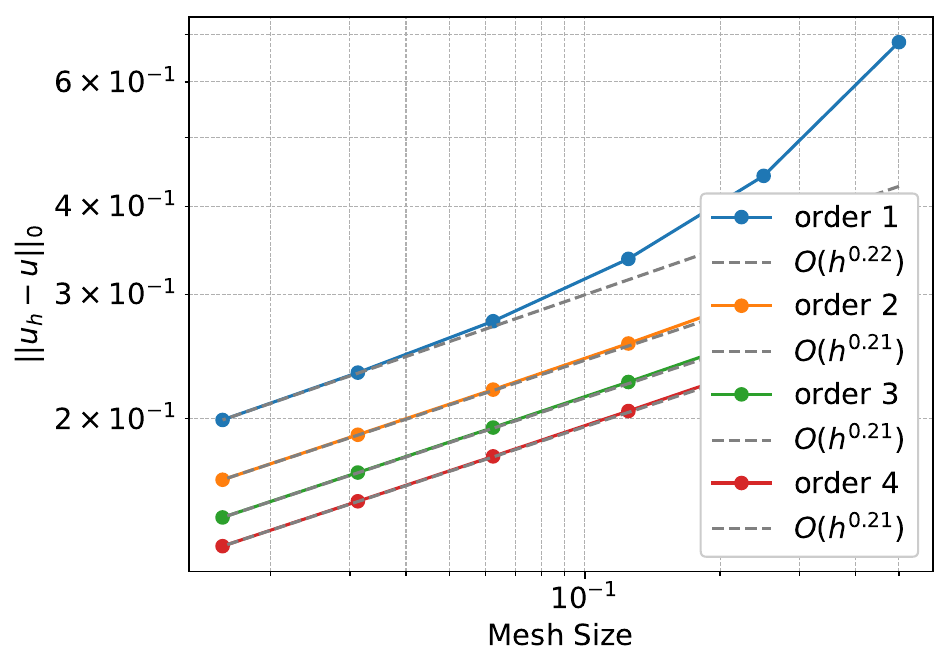}
  \includegraphics[width = 0.49\textwidth]{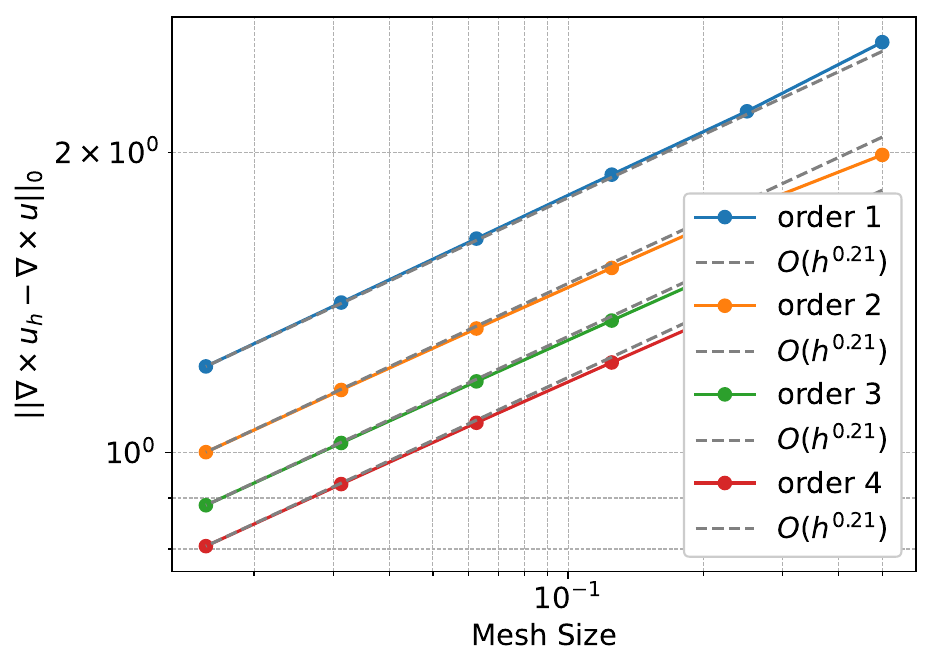} \\
  \includegraphics[width = 0.49\textwidth]{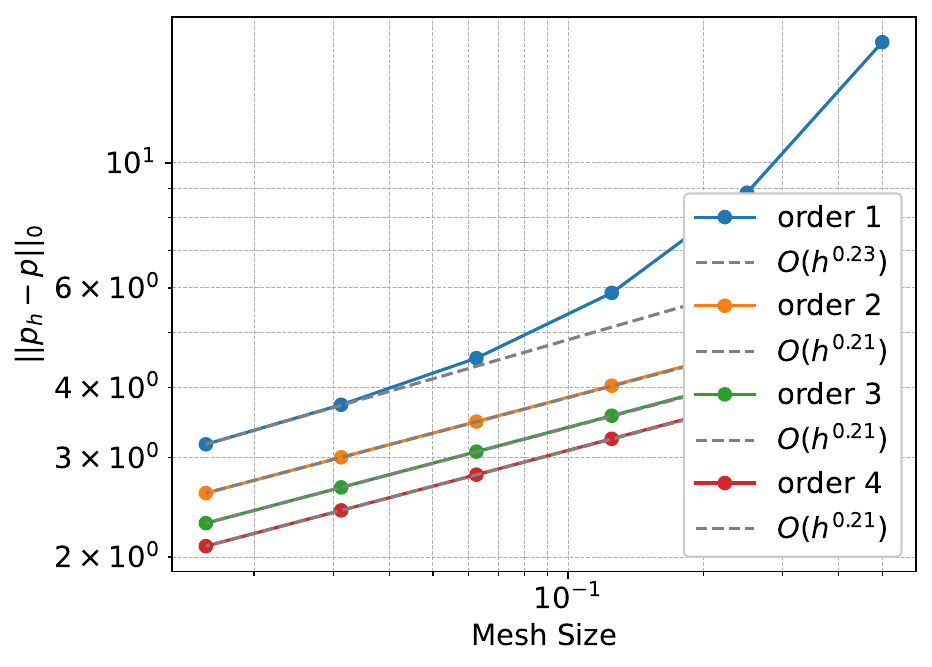}
  \caption{Convergence analysis for the non-stabilized scheme as described in \cref{sec:NonSmoothExperiment} of (top) $\bfu$ and $\nabla\times\bfu$ in the $L^2$ norm, and (bottom) $p$ in the $L^2$ norm for the experiment as discussed in \cref{sec:NonSmoothExperiment}. The results for lowest-order elements are labeled as \lq\lq order 1\rq\rq.}
  \label{fig:RobinConvergenceAnalysisLshapeNonStab}
\end{figure}

\begin{figure}[h]
  \centering
  \includegraphics[width = 0.49\textwidth]{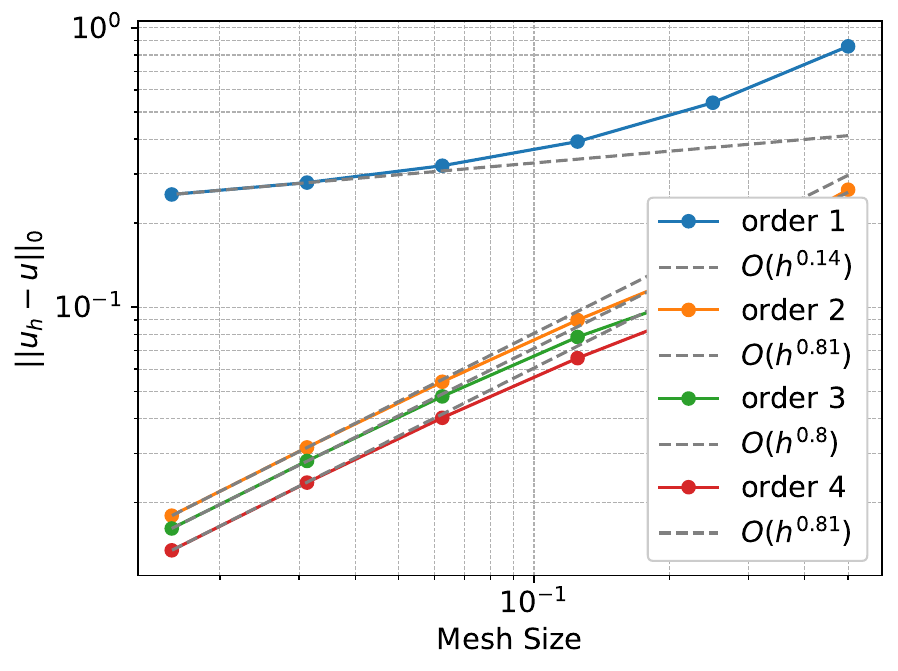}
  \includegraphics[width = 0.49\textwidth]{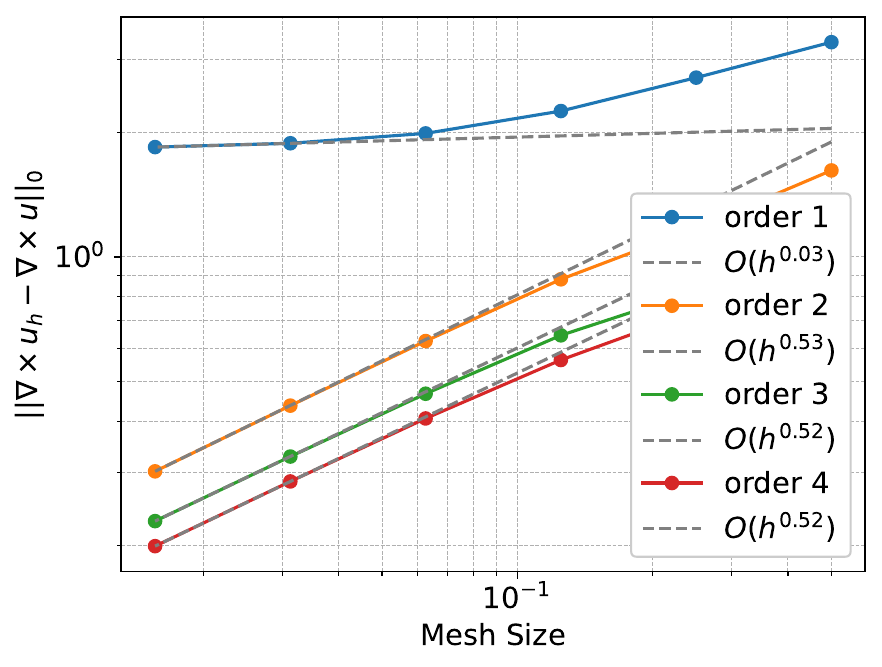} \\
  \includegraphics[width = 0.49\textwidth]{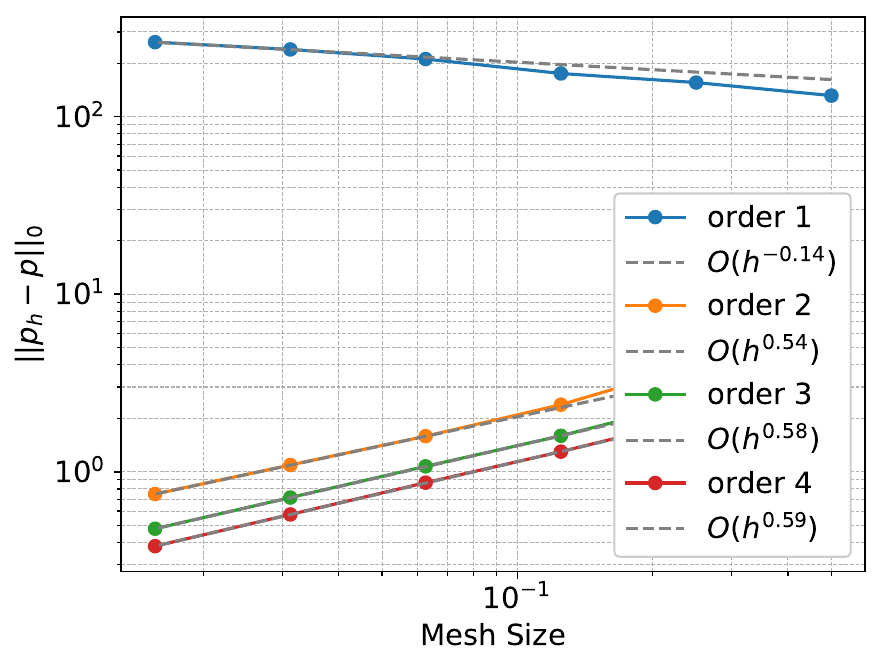}
  \caption{Convergence analysis for the stabilized scheme as described in \cref{sec:NonSmoothExperiment} of (top) $\bfu$ and $\nabla\times\bfu$ in the $L^2$ norm, and (bottom) $p$ in the $L^2$ norm for the experiment as discussed in \cref{sec:NonSmoothExperiment}. The results for lowest-order elements are labeled as \lq\lq order 1\rq\rq.}
  \label{fig:RobinConvergenceAnalysisLshapeStab}
\end{figure}

We display the $L^2$-error of $\bfu$, the $L^2$-error of $\nabla\times\bfu$, and the $L^2$-error of $p$ for the non-stabilized and stabilized scheme in \Cref{fig:RobinConvergenceAnalysisLshapeNonStab} and \Cref{fig:RobinConvergenceAnalysisLshapeStab}. While standard $H^1$-based formulations achieve convergence rates of approximately $\lambda \approx 0.54$ in the $H^1(\Omega)$ norm and $0.86$ in the $L^2(\Omega)$ norm (results not reported here for brevity), the $\Hcurl$ method attains a suboptimal rate of only about $0.21$ in all norms. We remark that this test falls outside the theoretical framework considered in this work, both because of the mixed slip--no-slip boundary conditions and because the domain is non-convex. At present, our developed theory is unable to explain the observed reduction in the convergence rate. Nevertheless, these results are compatible with the assumptions of \Cref{thm: optimal rates} being needed to guarantee the optimal convergence rates established therein, although several of them are violated simultaneously here. The stabilized version attains rates comparable to those of the $H^1$-based formulations, except for the lowest order. In fact, we observe a convergence rate of around 0.8 in the $L^2$-norm of the velocity and a convergence rate of around 0.54 in the other norms. For the lowest-order case, the stabilized method does not converge: the pressure error grows, while the velocity errors appear to stagnate, likely because the velocity space does not contain the full space of linear polynomials. This behavior is reminiscent of the results in \citep{WangYe}, where it is shown that lowest-order Raviart--Thomas elements lack sufficient approximation properties in the discrete $H^1(\Omega)$ norm.

The main takeaway of this section is that care is required when applying $\Hcurl$-based methods (e.g. those mentioned in \Cref{rem:non-convex}) to the Stokes problem on non-convex domains. Nevertheless, the addition of a simple stabilization term appears to restore both convergence to the correct solution and the convergence rates of the $H^1$-based formulations. However, a rigorous analysis of the stabilized formulation falls outside the scope of the present work and will be considered elsewhere.

\subsection{Lid-driven cavity} \label{sec:LidDrivenCavity}
\begin{figure}[htb]
  \centering
  \includegraphics[width = 0.48\textwidth]{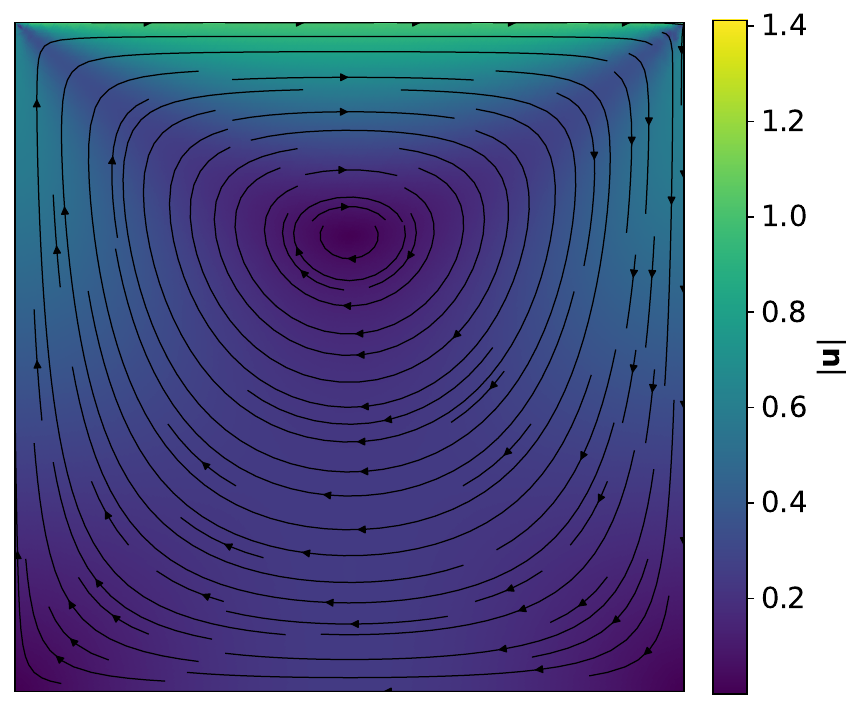}
  \includegraphics[width = 0.48\textwidth]{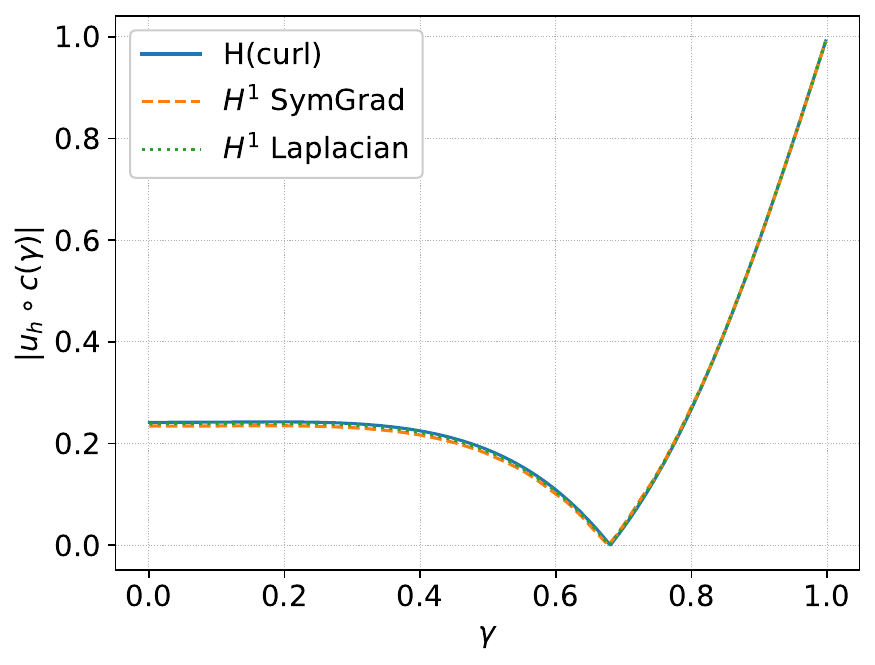}
  \caption{Visualization of the computed solution as described in \cref{sec:LidDrivenCavity} using 3rd order polynomials on an unstructured mesh with mesh-width $h=0.02$. The domain is a 1-by-1 square. (Left) The white lines represent the streamlines and the colors indicate the magnitude of the velocity field $\bfu$. (Right) The curves represent the magnitude of the discrete velocity $\bfu_h$ on the line segment $c(\gamma)=[0.5,\gamma]$, $\gamma\in[0,1]$.}
  \label{fig:LidDrivenCavity2D}
\end{figure}

As our third test case, we consider the lid-driven cavity. We follow \citep[Sect. 6.2, Test 2]{LimacheObjectivityTests} and define $\Omega$ as a unit square. We impose a Dirichlet boundary condition $\bfu = [1,0]$ on the top boundary (the \lq\lq lid\rq\rq) and Navier slip boundary conditions on the remaining sides. We use polynomials of degree $3$ on an unstructured mesh with mesh-width $h=0.02$. The results are presented in \cref{fig:LidDrivenCavity2D} (Left). Since all the boundaries are straight in this case, both the \lq\lq Divergence\rq\rq\, and the \lq\lq Laplace\rq\rq\, formulations converge to the same solution. In \Cref{fig:LidDrivenCavity2D} (Right) we plot the magnitude of the approximation of the velocity field $\bfu$ on the line segment $c(\gamma)=[0.5,\gamma]$. We see that the different computed solutions agree.

\subsection{Slippery annular cavity} \label{sec:Annulus}
\begin{figure}[htb]
  \centering
  \includegraphics[width = 0.48\textwidth]{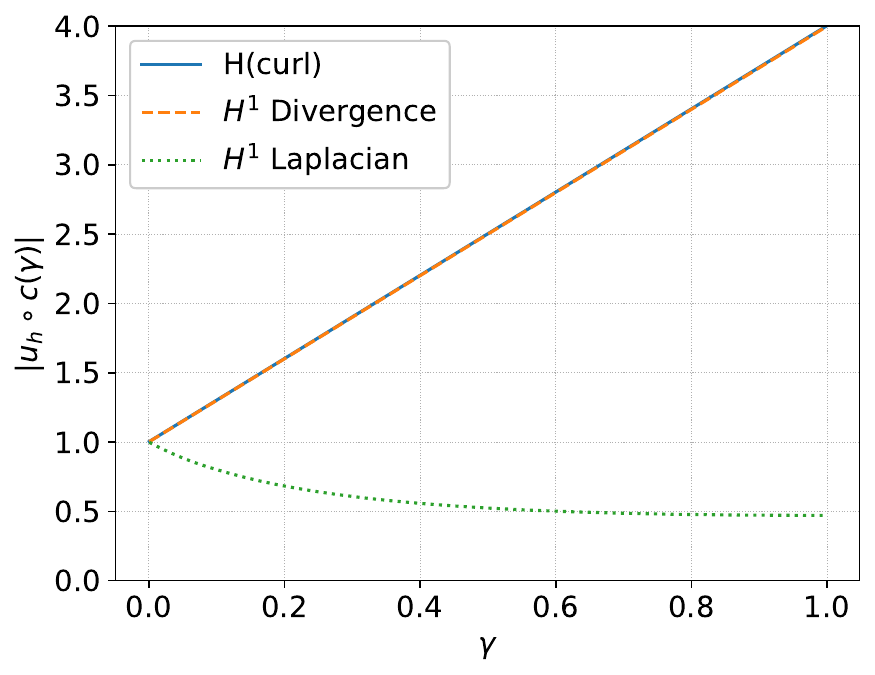}
  \includegraphics[width = 0.48\textwidth]{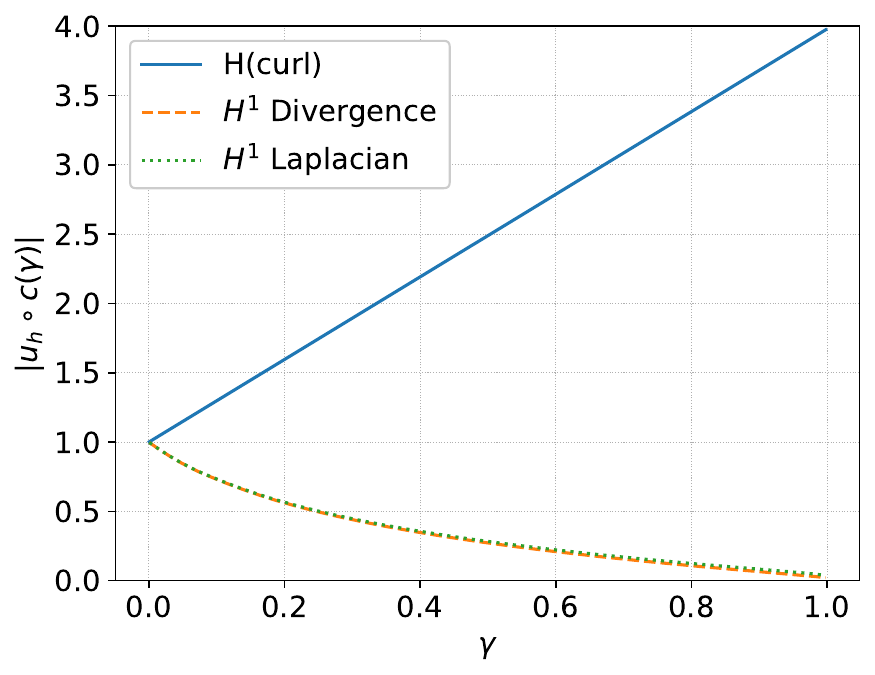}
  \caption{{Visualization of the computed solution along the line $c(\gamma)=[1+3\gamma,0]$, $\gamma\in[0,1]$, as described in \cref{sec:Annulus} on curved (Left) and non-curved (Right) meshes. The \lq\lq $\Hcurl$\rq\rq\, method correctly produces a rigid body motion on both curved and non-curved meshes. The lines represent the magnitude of the velocity field $\bfu_h$. Note that both the \lq\lq $\Hcurl$\rq\rq\, and \lq\lq Divergence\rq\rq\, produce the right solution on curved meshes, but only the \lq\lq$\Hcurl$\rq\rq\, produces the right solution on non-curved meshes. Lowest-order elements were used for all methods. }
  }
  \label{fig:Annulus2D}
\end{figure}

Next, we consider an annulus with inner radius of 1 and outer radius of 4 (see \cref{fig:Annulus2D}). On the inner boundary we impose Dirichlet boundary conditions $\bfu = [-y,x]$. On the outer boundary, we impose a Navier slip boundary condition. This test is of particular interest as it indicates whether the scheme satisfies the objectivity criteria from \citep[Sect. 6.5, Objectivity test 1]{LimacheObjectivityTests}. The objectivity criteria dictate that the true solution is a rigid rotation in this example. We used 1st order polynomials on an unstructured, curved mesh with mesh width $h=0.25$. As shown in  \Cref{fig:Annulus2D}, we observe that the numerical solution correctly represents a rigid body motion.

{Repeating the same experiment on the non-curved variant of the mesh, we observe on the right of \Cref{fig:Annulus2D} that our method still captures the correct rigid body motion, owing to the strategy proposed in \cite{NeunteufelCurvatureComputation}, which computes a distributional derivative of the normal vector. In contrast, in this case, the \lq\lq Divergence\rq\rq, formulation yields an incorrect solution on this mesh, owing to a Babu\v{s}ka-like paradox. The remedy, as shown in \citep{DioneEtAl, DioneEtAl2}, involves replacing the mesh normal vector with a regularized normal that accounts for the domain's curvature.}

\subsection{A manufactured solution in 3D} \label{sec:ManufacturedSolution3D}
\begin{figure}[htb]
  \centering
  \includegraphics[width = 0.49\textwidth]{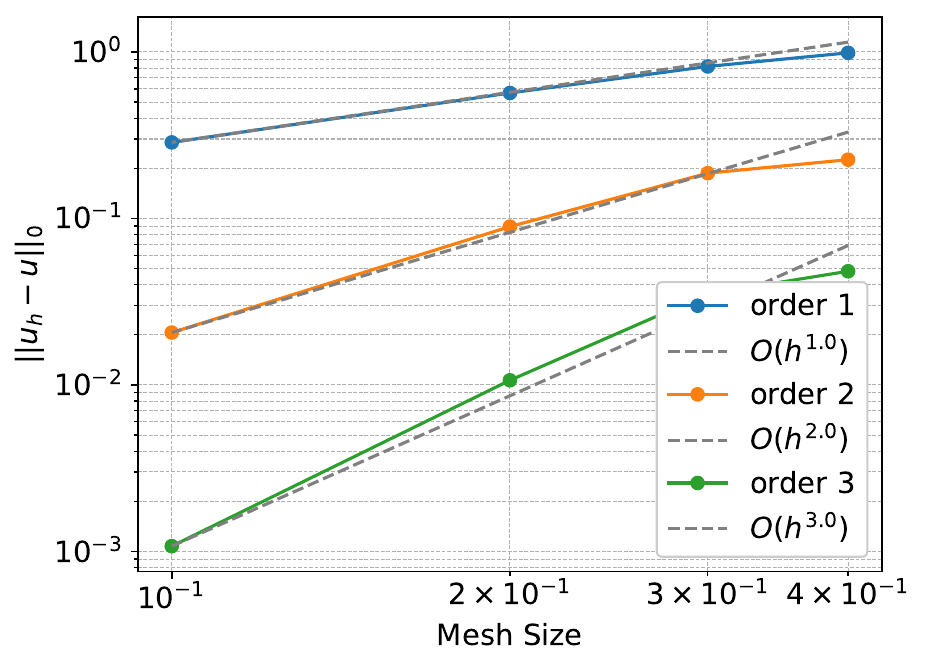}
  \includegraphics[width = 0.49\textwidth]{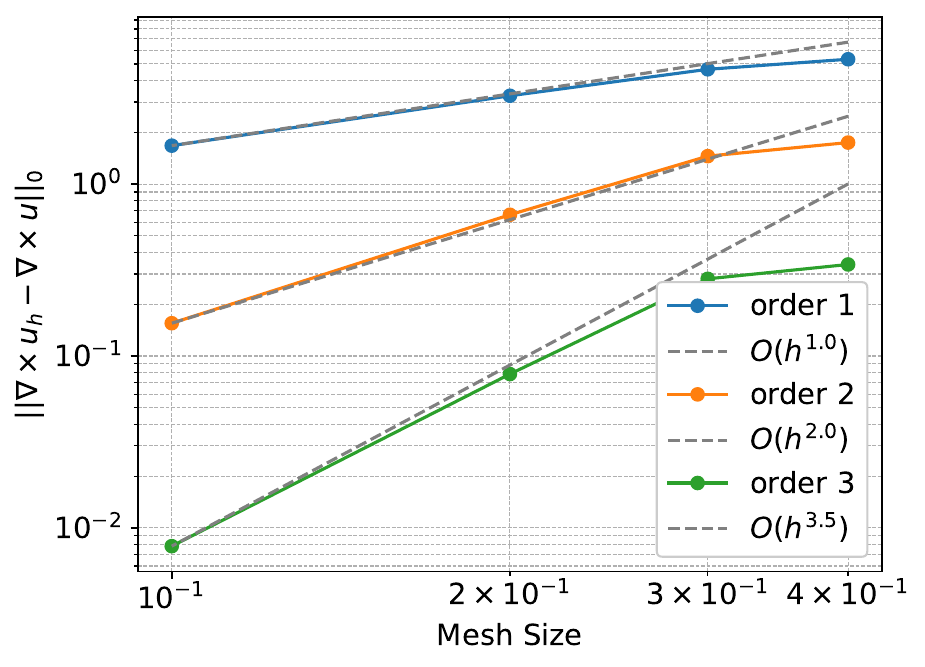} \\
  \includegraphics[width = 0.49\textwidth]{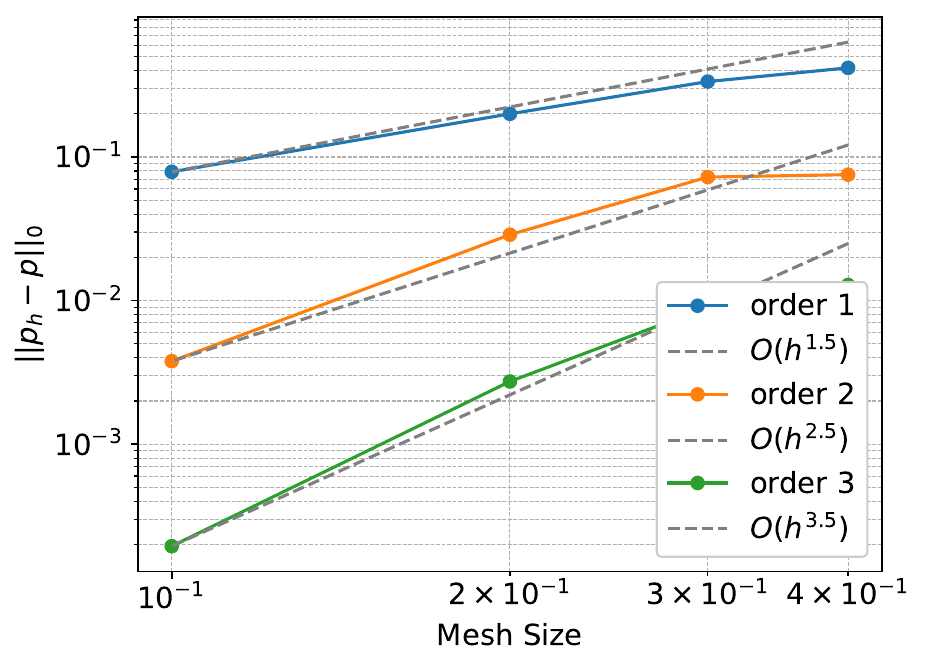}
  \includegraphics[width = 0.49\textwidth]{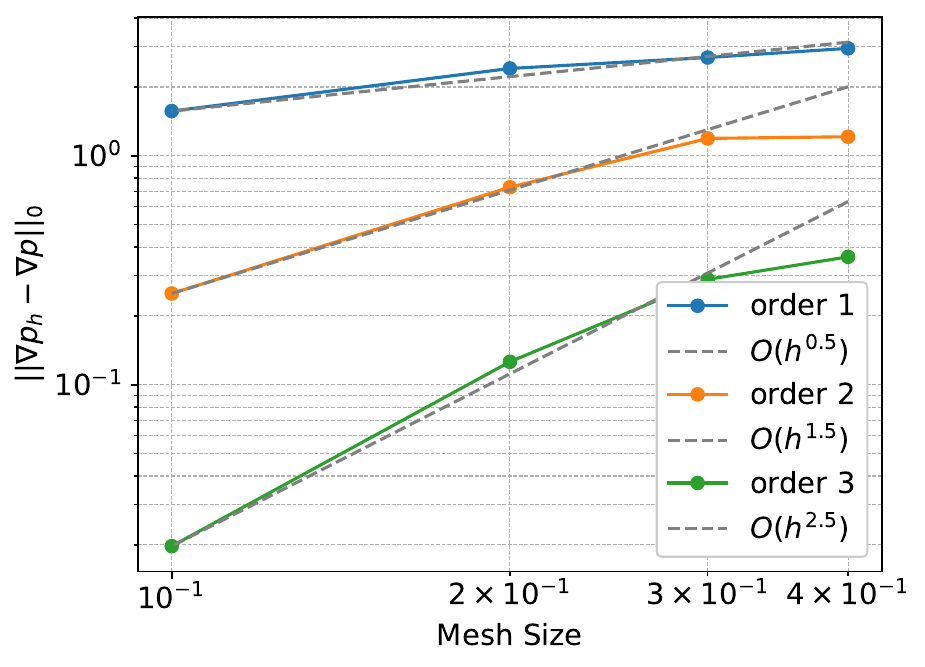}
  \caption{{Convergence analysis of (top) $\bfu$ and $\nabla\times\bfu$ in the $L^2$ norm, and (bottom) $p$ and $\nabla p$ in the $L^2$ norm for the experiment as discussed in \cref{sec:ManufacturedSolution3D}. The results for lowest-order elements are labeled as \lq\lq order 1\rq\rq.
  }}
  \label{fig:RobinConvergenceAnalysis3D}
\end{figure}

{In this experiment, we consider a manufactured solution on an ellipsoid $\Omega\subset\R^3$ with width $2$, height $1.6$, and depth $2.4$. Note that this domain has no continuous rotational symmetry and we can thus expect a unique solution. Given the solution
\begin{align*}
    \bfu &= \Curl\begin{bmatrix}
       \sin(\pi y)\cos(\pi z)x^2 \\
       \sin(\pi x)\cos(\pi x)yz\\
       \sin(\pi z)\cos(\pi x)xz^2
    \end{bmatrix}, &
    p = \frac{1}{10}\sin(\pi x)\cos(\pi y)\cos(\pi z),
\end{align*}
we derive $\mathbf{f}$, $z$, and $\mathbf{g}$ analogously to \cref{sec:experiment1}. Again, we approximate the Weingarten map as described in the beginning of this section. In \Cref{fig:RobinConvergenceAnalysis3D}, we display the $L^2$-error of $\bfu$, the $L^2$-error of $\nabla\times\bfu$, the $L^2$-error of $p$, and the $L^2$-error of $\nabla p$, respectively. The convergence rates match those predicted by \Cref{thm: optimal rates} in two dimensions, suggesting that the analysis could be extended also to three dimensions. As in the two-dimensional case, we observe half an order loss in convergence for the $L^2$ norm of $\nabla p$, consistent with the upper bound of \Cref{thm: error estimate p}.}
\subsection{Flow around a slippery sphere} \label{sec:FlowAroundSlipperySphere}
\begin{figure}[htp]
  \centering
  \includegraphics[width = 0.3\textwidth]{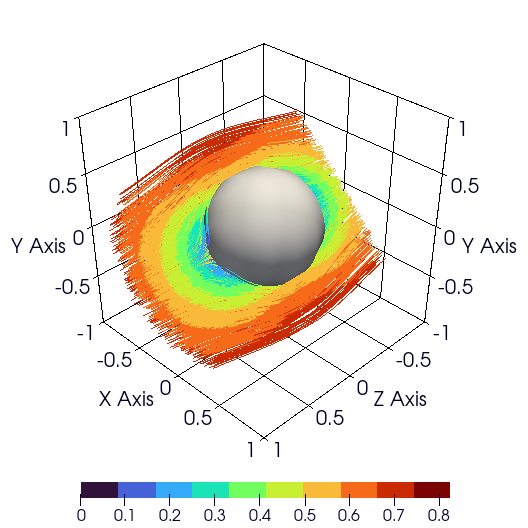}
  \includegraphics[width = 0.3\textwidth]{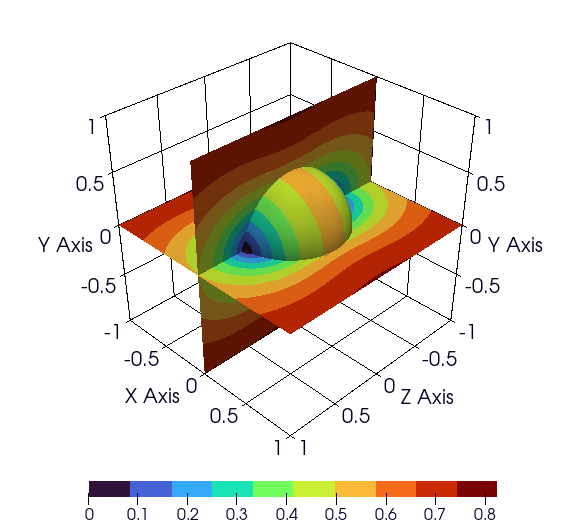}
  \includegraphics[width = 0.3\textwidth]{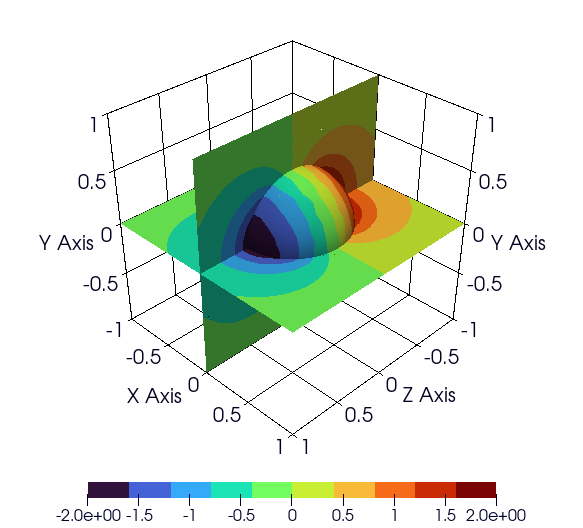}
  \caption{{Visualization of the computed solution as described in \cref{sec:FlowAroundSlipperySphere}. The domain is a 2-by-2-by-2 box with a sphere of radius 0.5 cut out at its center. We used 3rd order polynomials on an unstructured, curved mesh with mesh width $h=0.2$. (Left) The lines represent the streamlines and the colors indicate the magnitude of the velocity field $\bfu$. (Middle) The colors indicate the magnitude of the velocity field. (Right) The colors indicate the pressure.}}
  \label{fig:FlowAroundSlipperySphere}
\end{figure}

{
In this experiment, we reproduce a numerical experiment reported in \cite{CostaSlip}. We consider Stokes flow around a sphere in 3D. The domain is a 2-by-2-by-2 (bounding) box with a sphere of radius 0.5 at the center cut out from the box. On the boundary of the bounding box, we enforce Dirichlet boundary conditions, while on the boundary of the sphere we enforce homogeneous Navier slip boundary conditions. In spherical coordinates $(r,\theta,\varphi)$, with $\theta$ the polar angle measured from the $z$-axis, the exact solution of this problem is
\begin{equation*}
    u_r = \frac{ \cos(\theta)(r-0.5)}{r}, \quad u_{\theta} = \frac{ \sin(\theta)(0.5 - 2r)}{2r}, \quad u_{\varphi} = 0, \quad p = -\frac{ 0.5\cos(\theta)}{r^2}.
  \end{equation*}

We visualize the results in \cref{fig:FlowAroundSlipperySphere}. The computed velocity field agrees with \cite[Figure 15]{CostaSlip}. We also perform a convergence analysis and display the $L^2$-error of $\bfu$, the $L^2$-error of $\nabla\times\bfu$, the $L^2$-error of $p$, and the $L^2$-error of $\nabla p$. The convergence rates for $\bfu$ again match those predicted by \Cref{thm: optimal rates} in two dimensions. However, for $r\geq2$ we observe a full order loss in convergence for the $L^2$ norm of $\nabla p$. This is half an order worse than predicted by \Cref{thm: error estimate p}, while for $r=1$ the observed rate $1/2$ coincides with the predicted bound. This may be related to the mixed boundary conditions, which are outside our theory. For comparison, \cite{Nitsche} derives, for no-slip conditions imposed by Nitsche's method, an upper bound with a loss of three halves of an order. For $r\geq2$, the observed rate lies between the half-order loss of our upper bound and the three-halves-order loss of the upper bound in \cite{Nitsche}; explaining it would require a separate analysis. }

\begin{figure}[htb]
  \centering
  \includegraphics[width = 0.49\textwidth]{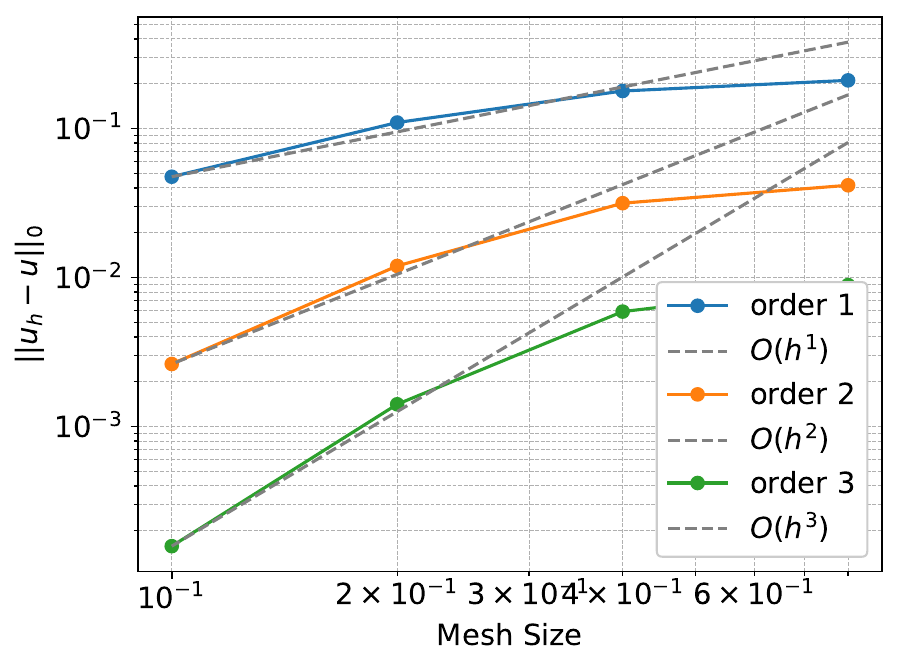}
  \includegraphics[width = 0.49\textwidth]{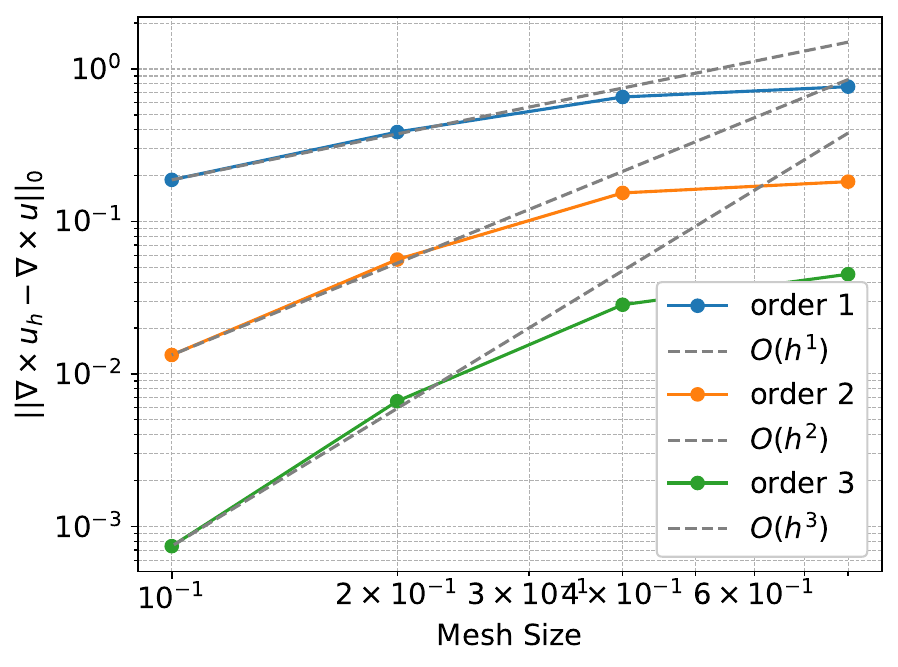} \\
  \includegraphics[width = 0.49\textwidth]{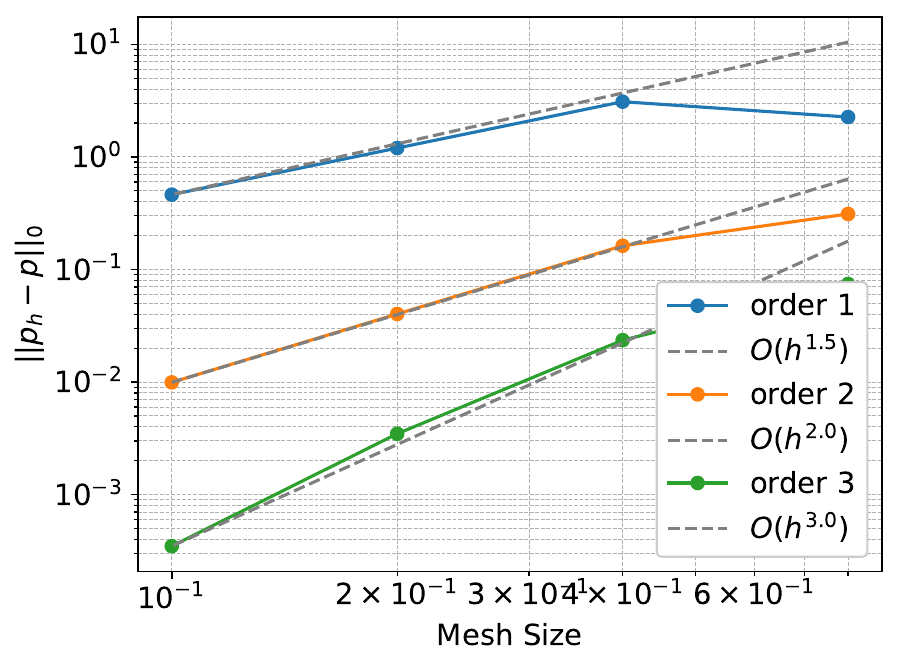}
  \includegraphics[width = 0.49\textwidth]{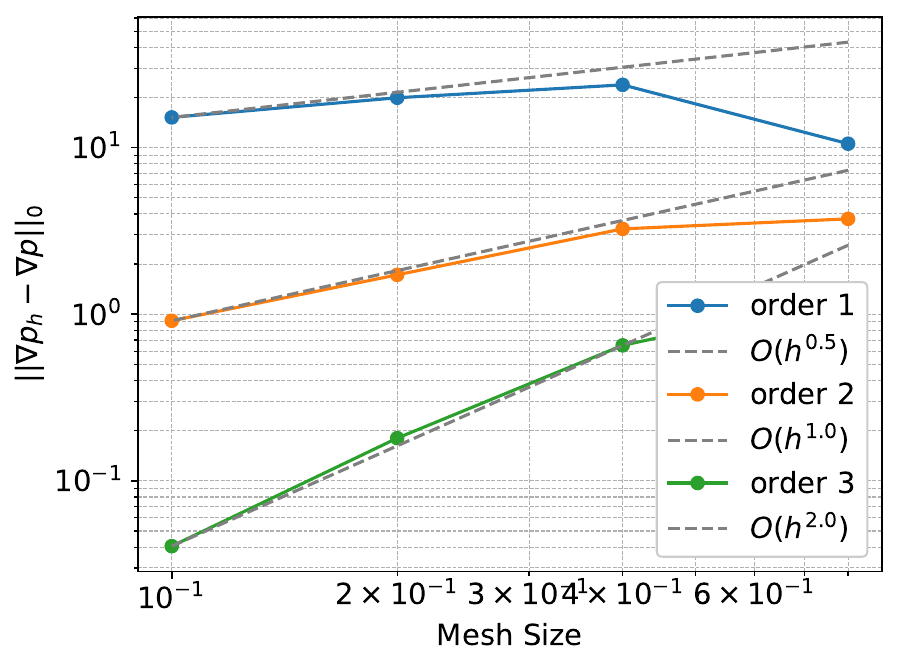}
  \caption{Convergence analysis of (top) $\bfu$ and $\nabla\times\bfu$ in the $L^2$ norm, and (bottom) $p$ and $\nabla p$ in the $L^2$ norm for the experiment as discussed in \cref{sec:FlowAroundSlipperySphere}. The results for lowest-order elements are labeled as \lq\lq order 1\rq\rq.
  }
  \label{fig:RobinConvergenceAnalysisSphere}
\end{figure}
\section{Conclusions}

We have treated the Navier slip boundary condition as a Robin condition for the Stokes
equations in rotational form. The associated Robin parameter is linked to the Weingarten
map of the boundary and may change sign along the boundary, so the resulting operator is
not necessarily semi-definite. Taking this feature into account, we proved that the
continuous problem admits a solution under a compatibility condition on the load (automatically satisfied when the kernel is trivial), unique up to the finite-dimensional kernel, with velocity in the subspace $\X\subset \Hcurl$
orthogonal to gradients.

At the discrete level, we considered a subspace of a
$\Hcurl$-conforming finite-element space that is only orthogonal to discrete gradients and
therefore not a subspace of $\X$. Consequently, the standard theory for conforming methods
did not apply. Instead we employed appropriate lifting and projection operators to obtain
stability and convergence estimates. These results were extended to an implementable
saddle-point formulation involving the pressure.

Numerical experiments confirmed the predicted convergence rates for a family of finite
element methods and further indicated that the discrete solutions exhibit physically
consistent behavior under the assumptions of the theory (convex or $\C^{1,1}$ domains); on a non-convex polygon, a stabilization appears necessary.

\section*{Acknowledgements}
This collaboration originated from discussions at the \emph{SuprenumPDE} workshop at the Bernoulli Center of EPFL in the summer of 2023.
The authors are grateful to Andrea Brugnoli for the stimulating discussions and for having
suggested the reference \citep{LimacheObjectivityTests}. E.Z. is a member of the GNCS-INdAM (Istituto Nazionale di Alta Matematica) group.

\section*{Declaration of generative AI use.}
ChatGPT (OpenAI, GPT-6) was used to assist with proofreading and improving the language and style of the manuscript. The authors take full responsibility for the content of the paper.

\appendix
\crefalias{section}{appendix}
\renewcommand{\theHsection}{appendix.\Alph{section}}
\renewcommand{\theHequation}{\theHsection.\arabic{equation}}
\section{A Nitsche-type method for Dirichlet boundary conditions} \label{sec: Nitsche}

Let $\Gamma = \Gamma_{\mathrm{D}} \cup \Gamma_{\mathrm{NS}}$, where $\Gamma_{\mathrm{D}}$ and $\Gamma_{\mathrm{NS}}$ denote the portions of the boundary on which Dirichlet and Navier slip boundary conditions are imposed, respectively. Following \citet{Nitsche}, we modify the bilinear form $a$ in \eqref{eq:continuous_problem} as follows:
\begin{align*}
    a(\bfu_h,\bfv_h) &\coloneqq
    (\nabla\times \bfu_h, \nabla \times \bfv_h)
    + \langle \alpha \ttrace(\bfu_h) ,\ttrace(\bfv_h)\rangle_{\Gamma_{\mathrm{NS}}} \\
    &+ \langle \nvec \times \nabla\times \bfu_h, \bfv_h \rangle_{\Gamma_{\mathrm{D}}}
    + \langle \bfu_h, \nvec\times \nabla \times \bfv_h\rangle_{\Gamma_{\mathrm{D}}}
    + \frac{C_w}{h} \langle \nvec \times\bfu_h, \nvec \times \bfv_h\rangle_{\Gamma_{\mathrm{D}}}.
\end{align*}
Here, $C_w>0$ is a penalty parameter, which has to be chosen sufficiently large, see \citet{Nitsche}.

\section{Calculated example on a cube}
\label{sec:calculated_example_on_a_cube}

Let $\Omega = (0,1)^3$ and define $\mathbf{A} : \Omega \to \mathbb{R}^3$ and $\bfu = \Curl \mathbf{A}$ as
\begin{align}
    \mathbf{A} &= \begin{pmatrix} 0 \\ \sin(\pi x)\sin(\pi z) \\ \sin(\pi x)\sin(\pi y) \end{pmatrix}, &
    \bfu &= \pi\begin{pmatrix} \sin(\pi x)[\cos(\pi y) - \cos(\pi z)] \\ -\cos(\pi x)\sin(\pi y) \\ \cos(\pi x)\sin(\pi z) \end{pmatrix}.
\end{align}

Since $\operatorname{div} \bfu = 0$ and its normal trace vanishes on each face of $\Gamma$, it follows that $\bfu\in \Hndivz$. Moreover, the symmetric gradient of $\bfu$ is given by
\begin{align}
    \symgrad(\bfu) &= \pi^2 \cos(\pi x)
    \begin{pmatrix}
        \cos(\pi y) - \cos(\pi z) & 0 & 0 \\
        0 & -\cos(\pi y) & 0 \\
        0 & 0 & \cos(\pi z) \end{pmatrix},
\end{align}
from which we conclude that $[\symgrad(\bfu)\nvec]_{\mathrm{tan}} = 0$ on $\Gamma$. Finally, we compute
\begin{align}
    - 2 \Div \symgrad(\bfu) =
    {2} \pi^3 \begin{pmatrix} \sin(\pi x)[\cos(\pi y) - \cos(\pi z)] \\ -\cos(\pi x)\sin(\pi y) \\ \cos(\pi x)\sin(\pi z) \end{pmatrix}
    = 2 \pi^2 \bfu.
\end{align}

From these properties, we conclude that $\bfu$ is the solution to \eqref{eq:continuousSymGradStokes} with $p = 0$ and $\mathbf{f} = 2 \pi^2 \bfu$. In turn, \Cref{thm:WvsW} implies that $\bfu$ is the solution to \eqref{eq:continuous_problem} and \Cref{rem:flat faces,rem:corners} further state that this holds with $\alpha = 0$.
We validate these claims by choosing $\bfv = \bfu$ in \eqref{eq:continuous_problem} and computing as follows:
\begin{align}
    \langle \alpha \ttrace(\bfu), \ttrace(\bfu)\rangle
    = (\mathbf{f}, \bfu) - (\Curl \bfu, \Curl \bfu) \
    &= 2 \pi^2 \| \bfu \|_{\VectorLtwo}^2 - \| \Curl \bfu \|_{\VectorLtwo}^2 \nonumber \\
    &= 2 \pi^4 - 2 \pi^4
    = 0{.}
\end{align}
Since $\ttrace(\bfu)$ is non-zero almost everywhere on $\Gamma$, we conclude that the energy identity $a(\bfu, \bfu) = (\mathbf{f}, \bfu)$ is satisfied if $\alpha$ is zero. {This example shows that we do not need to introduce a singular Dirac term on the edges of the cube. On the edge aligned with the $z$-axis, for example, the boundary conditions already ensure that $\bfu \cdot \mathbf{e}_1 = \bfu \cdot \mathbf{e}_2 = 0$, which annihilates any term concerning those components. The energy identity holds with $\alpha=0$ and no additional edge term; a singular term with nonzero energy for $\bfu \cdot \mathbf{e}_3$ would violate it.}

\end{document}